%% file: SIAPLS-2022.tex
\documentclass[final,12pt]{article}
\input{PrePrintMacros.tex}

\begin{document}

\maketitle

\begin{abstract}
\input{Manuscript/abstract.tex}
\end{abstract}

{\bf Keywords}: \input{Manuscript/keywords.tex}
\vspace{.25cm}

{\bf MSC}: 
\input{Manuscript/ams.tex}

\graphicspath{{Figures/}}

\input{Manuscript/intro.tex}

\input{Manuscript/mechanics.tex}

\input{Manuscript/discretization.tex}

\input{Manuscript/strainenergies.tex}

\input{Manuscript/experiments.tex}

\input{Manuscript/conclusions.tex}

\input{Manuscript/acknowledgments.tex}

\input{Manuscript/note.tex}

\bibliography{references.bib}

\input{Manuscript/SIAMLS_Supplement}
\end{document}

%% file: PrePrintMacros.tex
\usepackage[utf8]{inputenc}
\usepackage{anysize} 
\marginsize{2cm}{2cm}{2cm}{2cm} 
\usepackage{fancyhdr} 
\usepackage{lipsum} 
\usepackage{color} 
\usepackage{bm,bbm} 
\usepackage{amsfonts,amsthm,amsmath,amssymb} 
\usepackage{latexsym,graphicx,epstopdf} 
\usepackage{etoolbox} 
\usepackage{hyperref} 
\usepackage[left,mathlines]{lineno} 
\renewcommand{\a}{{\bf a}}
\renewcommand{\b}{{\bf b}}

\renewcommand{\d}{{\bf d}}
\newcommand{\e}{{\rm e}}
\newcommand{\f}{{\bf f}}

\newcommand{\n}{{\bf n}}

\renewcommand{\u}{{\bf u}}
\renewcommand{\v}{{\bf v}}
\newcommand{\w}{{\bf w}}
\newcommand{\x}{{\bf x}}

\newcommand{\zero}{{\bf 0}}

\newcommand{\cc}{{\bf C}}

\newcommand{\ff}{{\bf F}}

\newcommand{\hh}{{\bf H}}
\newcommand{\ii}{{\bf I}}

\newcommand{\pp}{{\bf P}}
\newcommand{\qq}{{\bf Q}}

\renewcommand{\ss}{{\bf S}}

\newcommand{\uu}{{\bf U}}
\newcommand{\vv}{{\bf V}}

\newcommand{\xx}{{\bf X}}

\newcommand{\iii}{\mathbb{I}}

\newcommand{\ppp}{\mathbb{P}}
\newcommand{\ppm}{\mathbbm{p}}

\newcommand{\rrr}{\mathbb{R}}

\newcommand{\tttt}{\mathcal{T}}

\newcommand{\btau}{{\bm\tau}}
\newcommand{\bxi}{{\bm\xi}}
\newcommand{\bchi}{{\bm\chi}}
\newcommand{\bsigma}{{\bm\sigma}}

\newcommand{\bDiv}{{\bf Div}}

\newcommand{\bdiv}{{\bf div}}

\newcommand{\tr}{{\rm tr}}
\newcommand{\transpose}{\texttt{t}}

\newcommand{\pder}[2]{\dfrac{\partial #1}{\partial #2}}
\newcommand{\Dder}[2]{\dfrac{\mathrm{D} #1}{\mathrm{D} #2}}

\newcommand{\iprod}[1]{\left( #1 \right)_{V_0}}

\newcommand{\myemail}{s.dominguez@usask.ca}

\newcommand{\funding}{This work was partially supported by CONICYT-Chile, the Pacific Institute for the Mathematical Sciences, and NSERC-Canada through the Discovery program.}

\title{A three-dimensional model of skeletal muscle tissues\thanks{\funding}}

\author{Javier A. Almonacid\footnotemark[3]
\quad 
Sebasti\'an A. Dom\'inguez-Rivera\thanks{Department of Mathematics and Statistics, University of Saskatchewan, Canada (\href{mailto:\myemail}{\myemail}) }
\quad 
Ryan N. Konno\footnotemark[3]
\\
Nilima Nigam\thanks{Department of Mathematics, Simon Fraser University, Canada 
(\href{mailto:javiera@sfu.ca}{javiera@sfu.ca}{javiera@sfu.ca}, \href{mailto:rkonno@sfu.ca}{rkonno@sfu.ca}, \href{mailto:nigam@math.sfu.ca}{nigam@math.sfu.ca}, \href{mailto:cassidyt@sfu.ca}{cassidyt@sfu.ca}).}{\,\,\,\thanks{Corresponding author: 
\href{mailto:nigam@math.sfu.ca}{nigam@math.sfu.ca}.}}
\quad 
Stephanie A. Ross \thanks{Department of Physical Therapy, University of British Columbia, Canada (\href{mailto:stephanie.ross@ubc.ca}{stephanie.ross@ubc.ca})}
\quad 
Cassidy Tam\footnotemark[3]
\quad 
James M. Wakeling \thanks{Department of Biomedical Physiology and Kinesiology, Simon Fraser University, Canada 
(\href{mailto:wakeling@sfu.ca}{wakeling@sfu.ca})}
}

%% file: Manuscript/abstract.tex
Skeletal muscles are living tissues that can undergo large deformations in short periods of time and that can be activated to produce force. In this paper we use the principles of continuum mechanics to propose a dynamic, fully non-linear, and three-dimensional model to describe the deformation of these tissues. We model muscles as a fibre-reinforced composite and transversely isotropic material. 
We introduce a flexible computational framework to approximate the deformations of skeletal muscle to provide new insights into the underlying mechanics of these tissues. The model parameters and mechanical properties are obtained through experimental data and can be specified locally. A semi-implicit in time, conforming finite element in space scheme is used to approximate the solutions to the governing nonlinear dynamic model. We provide a series of numerical experiments demonstrating the application of this framework to relevant problems in biomechanics, and also discuss questions around model validation.


%% file: Manuscript/keywords.tex
skeletal muscle, non-linear deformation, continuum mechanics, hyperelasticity, dynamic contraction, isometric contraction, finite element method, time stepping method

%% file: Manuscript/ams.tex
74B20, 74L15, 65M60

%% file: Manuscript/intro.tex
\section{Introduction}\label{section:3dmuscleintro}
Why do animals move the way they do? Do the muscles of mice and giraffes and frogs behave identically? Is the force produced by a muscle (scaled by size) uniform across species?  What transpires during neuromuscular disorders, and of the many changes which we see \--- which are the most significant? Such questions continue to motivate large numbers of researchers, across the fields of biology, physiology, biomechanics, materials science, and medicine. A full understanding of the mechanical role of muscle architecture and the properties  of its tissues remains an open and challenging scientific problem. This article presents some recent work involving the mathematical modeling and simulation of skeletal muscle mechanics.  Of course, given the complexity of animal locomotion and muscle behavior, numerous important questions are not addressed in this work. We hope this article sparks the interest of readers and encourages the involvement of more mathematical research into this fascinating scientific area. 

Skeletal muscles are highly complex tissues which play a crucial structural and mechanical role in living beings.  Unlike non-living materials, muscles can be deformed both {\it actively} or {\it passively}. Active deformations are experienced when internal electro-chemical processes and
changes in the shape of proteins are triggered by neuronal signals (see, e.g. \cite{ref:oomens2009,ref:wakeling2020}). Passive deformations describe the non-linear mechanical response of muscle tissue in the presence of external forces \cite{ref:hill1938}. Both voluntary and involuntary movements are possible within skeletal muscle; the cells within skeletal muscle are long and multi-nucleated, and respond to external neuronal stimuli. In contrast,  cardiac muscle (which has also been extensively studied within the mathematics community) consists of mono-nucleated cells which contract on their own intrinsic rhythm; see, for instance \cite{Mescher2018}.

Given the complexity of skeletal muscle, teasing apart the role played by shape (the length, cross-sectional area, angle of attachment to joints), architecture (the arrangement of constitutive tissues), pennation (the angle at which muscle fibres are arranged), the arrangement of fast and slow-twitch fibres, which fibres are enervated and the mechanical properties of the individual tissues (including their their stiffness, compressibility, density) is notoriously challenging in experiments. Several of these change as we age \--- there is increased fat infiltration in our muscle, we lose muscle mass, connective tissues become stiffer. Which of these contributes more significantly to loss of muscle force output, or are these all equally important effects? Likewise, neuromuscular diseases lead to changes in architecture, shape, muscle properties and overall muscle mass. It is near-impossible to vary these properties within a given test subject experimentally to determine the differential impact of these changes. The goal of therapeutic treatments is to address the most significant factors impacting muscle function. At the same time, advances in medical imaging provide us with the possibility of generating patient-specific computational meshes on which we can examine the role of these factors. The framework we present in this paper has been used to study some of these questions, and in this paper we describe this in detail. One of the studies included here is on integrating MRI-based data to study muscle contraction; another is on furthering our understanding of cerebral palsy. We believe that advances in experimental science, medical imaging and computational power can be powerfully used in tandem to further our understanding of muscle mechanics.

One-dimensional models of skeletal muscle have been developed to achieve a deeper understanding of skeletal muscle (e.g. \cite{ref:hill1938,ref:zajac1989}). These models, generally referred to as {\it Hill-type muscle models}, have provided a substantial insight into the overall behavior of skeletal muscle. This includes insight in the basics of muscle mechanics \cite{ref:zajac1989}, the role of muscle in the musculoskeletal system \cite{Delp2007}, and the impact of muscular diseases on force output \cite{VanDerKrogt2016}. While these models are simple and reasonably accurate in some situations (particularly when considered as part of a full musculoskeletal study), they fail to account for many features of whole muscle dynamics and architecture. Mass and inertial effects are not considered in most Hill-type models, and indeed may be negligible if considering single-fibre experiments. It is frequently assumed that larger muscles behave as scaled-up single fibres (see, e.g. \cite{ref:zajac1989}). However, in larger muscle the impact of mass cannot be neglected, especially when the fibres are sub-maximally activated  (see eg. \cite{ref:holt2014,ref:ross2016,ross2020}). In addition, the three-dimensional architecture of muscle cannot be represented by Hill's model alone. Muscles are nearly incompressible tissues \cite{ref:baskin1967}, and thus, when muscle deforms, parts of the surface of the muscle deform inwards while other parts deform outwards, even during isometric contractions \cite{ref:randhawa2018,ref:ryan2020,ref:wakeling2020}. In this paper an {\it isometric contraction} is defined as  an elastic deformation of the tissues in which two opposing ends of the muscle are fully clamped. In the physiology literature however, isometric contraction refers to the combined length of the muscle and tendon being invariant. Additionally, the deformation also depends on other tissues not explicitly considered in many one-dimensional models; aponeurosis, base material (comprising extracellular materials and fluids), and fat all affect the dynamics and mechanics of muscles, see e.g. \cite{ref:rahemi2015,ref:ross2018-2,ref:bollinger}. Therefore, to be able to better capture and understand the behavior of muscle deformation in three dimensions, a more involved model is needed. 

To address the aforementioned concerns and fully capture the mechanics behind skeletal muscle, three-dimensional models have been developed. Many of these models are phenomenological continuum models that are developed using whole muscle properties (e.g. \cite{ref:blemker2005,Hedenstierna2008,Heidlauf2014,Heidlauf2017,Knaus2022,Li2019,Li2021A,Maier2022, ref:oomens2003,Osth2012,rahemi2014regionalizing,Rohrle2012,Roux2021,ref:wakeling2020,Webb2014,Weickenmeier2014,ref:yucesoy2002}). We note in particular \cite{ref:blemker2005} (a quasi-static approach), which has lead to an explosion of three-dimensional muscle models coupled with finite element simulations. These models are particularly useful as they are able to capture the architecture and bulging of muscle. Furthermore, these models have been developed in both quasi-static and dynamic formulations. While these models have been able to accurately capture the overall mechanics of skeletal muscle, they are computationally expensive. This makes them difficult to implement in larger musculoskeletal models such as OpenSim \cite{Delp2007}. Recently, authors in \cite{Gfrerer2021} developed a fibre based model, which has a lower computational cost than the typical continuum models, for the purposes of implementing more complex simulations with multiple muscles.
While these models focus on the mechanics of the whole muscle, to reach a deeper understanding of how muscle functions, the influence of the microstructure needs to be considered. 

The length scales involved in muscle mechanics range from the $\mu$m for sarcomeres  to the size of the gastrocnemius in larger mammals (in the order of meters) it is computationally expensive to model the behavior of large muscles by considering the individual behaviors of each filament. Micromechanical models have been developed to capture the influence of features, such as the connections of the cross-bridges, on the macroscopic scale \cite{ref:bol2008,Lamsfuss2021,ref:rohrle2007}.Another method of capturing microscopic mechanics on a macroscopic scale is the use of homogenization methods (e.g. \cite{Bleiler2019,Konno2021,ref:spyrou2017}). These models are utilized to capture the influence of the microstructure on the macroscopic scale in a computationally efficient way. Additionally, some continuum models have been developed to look at force transmission from microscopic features to the macroscopic structures \cite{ref:sharafi2010}.

Skeletal muscle models have also been developed to specifically address the mechanics of impaired muscle tissue, such as that affected by disease, disuse, or aging. For example, authors in \cite{Stefanati2020} developed a chemo-mechanical model of muscle tissue that captures the influence of muscular dystrophy, and the implications of the disease on the fibre level have been investigated in \cite{Virgilio2015}. Continuum models have yielded qualitatively similar results around regional deformations as experimental data of maximal contractions at constant muscle velocity, as well as during fixed-end contractions (e.g. \cite{ref:ross2018-2}). Further, these models can be used to understand the influence of transverse compression \cite{ref:ryan2020} and fat infiltration \cite{ref:rahemi2015} on muscle force.

The continuum models of muscle mentioned before focus on quasi-static deformation, and do not include any velocity-dependent effects; however, the velocity-dependent creep and relaxation have been found to be important in muscle mechanics \cite{ref:hill1938}. Viscoelastic models of muscle have been developed to capture these effects. \cite{ref:vanloocke2008} developed a viscoelastic model to look at the time-dependent properties of muscle in passive muscle tissue. More recently, \cite{Ahamed2016} developed a passive viscoelastic model using experimental data for the cyclic lengthening tests. These models do not, however, look at the active mechanics of muscle tissue using a viscoelastic model. Other dynamic effects that need to be considered are the effects of mass \cite{ref:ross2016,ref:ross2018-1,ref:ross2018-2}. Muscle deformation is inherently dynamic, and so in this paper we describe a fully dynamic model of skeletal muscle that will capture the velocity-dependent contributions to muscle mechanics in three dimensions. Our approach is  focussed on the meso-scale. With this model we are able to capture the influence of different material types found in the {\it muscle-tendon unit} (MTU): muscle, aponeurosis, fat, and tendon. We have performed computational investigations involving the first three of these both in physiological and pathological conditions (see e.g. \cite{ref:ross2018,ref:rahemi2015}), while the inclusion of tendon is part of our future work. Each of these materials possess different mechanical and physiological properties, and understanding the influence of the individual materials is critical when considering diseased muscle tissue \cite{Konno2022}.

Most of the current knowledge on muscle mechanics comes from experiments performed in laboratories, and generally on mammals. An important part of our work is to be able to fit the models that describe the mechanics of each of the tissues in skeletal muscle to experimental data, generally available in the form of strain-stress  data points. These data are collected under different experimental conditions, and from different species. As a consequence, we restricted our attention to parameters in the model which were significant, and have tested the sensitivity of computations to variations in them. 
Validation of the overall model is key, particularly to avoid over-fitting. This is later discussed in more depth in the paper.

In order to properly capture many different biomechanical phenomena, we designed our numerical framework to be modifiable to a variety of different experimental settings. An important design principle was to accurately capture the dynamics, while still being able to model a quasi-static simulation of muscle contraction. \textit{in silica} experiments have been performed with the use of our framework for isometric contractions to study the energetics of muscle contractions \cite{ref:wakeling2020}, and of compression loads on muscles \cite{ref:ryan2020}. In the dynamic case, the framework has been used to simulate cyclic contraction \cite{ref:Ross2021,ref:ross2018-2}.
With this framework it is possible to specify different mechanical properties of various tissues (muscle, aponeurosis,  fat, base material and tendon), localized fibre orientation vectors, regionalized activation of muscle fibres, as well as importing custom-made (MRI-derived) meshes.

The rest of this paper is organized as follows. In \autoref{section:3dmusclemechanics} we provide an overview of muscle physiology, non-linear elasticity, and introduce the system of equations that we use to describe the large deformation of muscles. We also discuss the role of muscle mass. 
In \autoref{section:3dmusclemethods} we describe the methods used to determine  the discretization strategy used for the model equations. A discussion of parameter fitting to data is presented in \autoref{sec:ParameterFitting}.  Finally, some results to validate our framework and demonstrate its biomechanical and clinical significance are presented in \autoref{section:3dmuscleexps}.


%% file: Manuscript/mechanics.tex
\section{Muscle physiology, mechanics, and modelling}\label{section:3dmusclemechanics}
In this section we first overview skeletal muscle physiology and include a description of the classical Hill-type muscle model \cite{ref:hill1938}. Then we formulate a three-dimensional dynamic skeletal muscle model, as well as its quasi-static counterpart.

\subsection{Skeletal muscle physiology}
While skeletal muscle structure varies from micrometers to the scale of meters, the modelling approach in this study will be focused on a macroscopic description of muscle and considers {\it muscle fibres} as the smallest unit. Muscle fibres are long, cylindrical cells composed of contractile unites called \textit{sarcomeres} arranged in parallel; the sarcomere is the functional unit of muscle fibres, which contract to produce force. The sarcomeres are on the scale of micrometers (typically about $2.2 \mu$m in length \cite{burkholder2001sarcomere}). They consist of thick myosin filaments interdigitated with thinner actin filaments, which are connected together via cross-bridges. The force production by the sarcomere can be characterized using the cross-bridge theory developed by A.F Huxley \cite{Huxley1957} and H.E. Huxley\cite{HuxleyHE2000}. Between each sarcomere is the Z-disk, which connects to the myosin through a protein called titin. These proteins are largely responsible for any passive response from the muscle fibres (see e.g. \cite{Herzog2018}).

\begin{figure}[!ht]
    \centering
    \includegraphics[width = 0.6\textwidth]{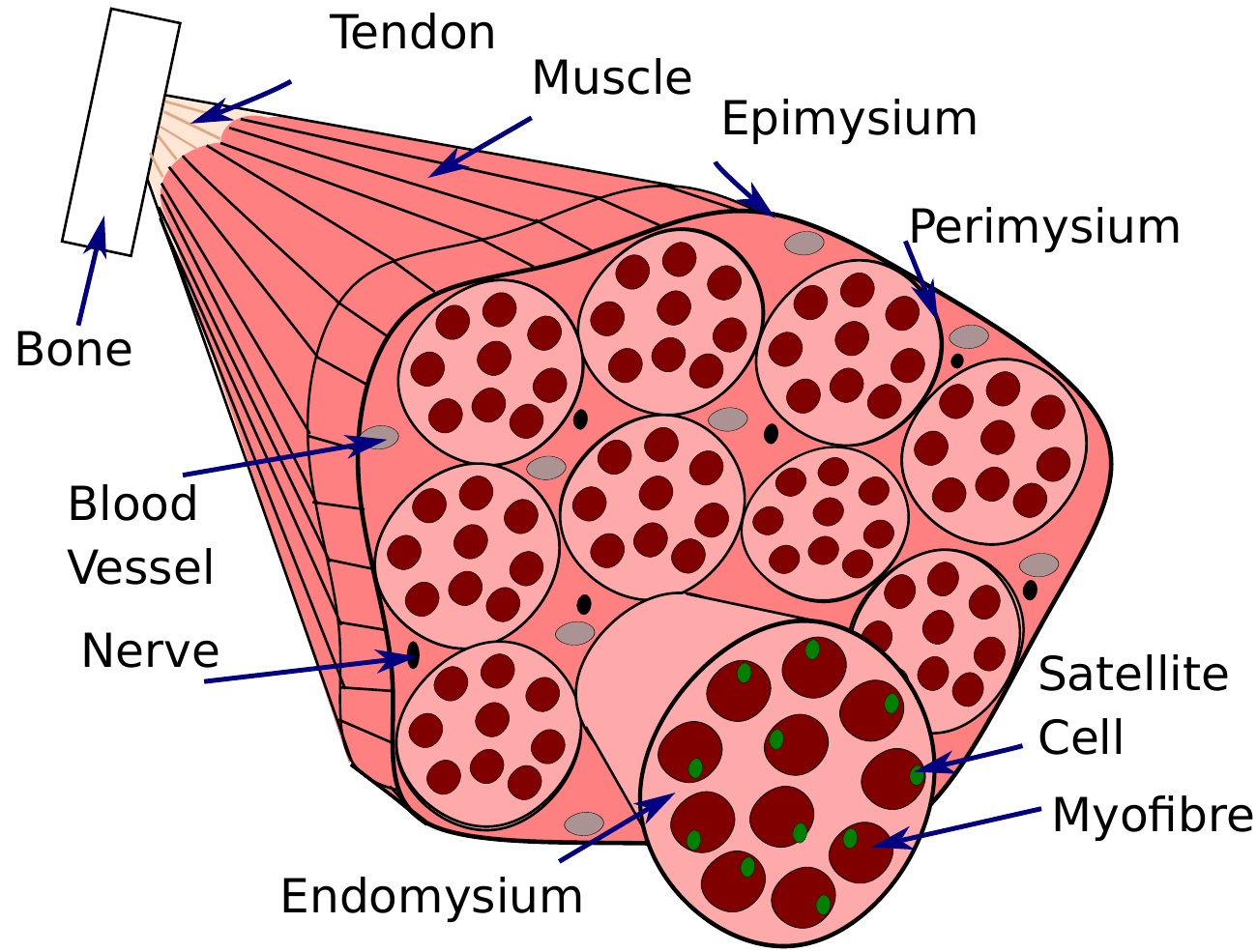}
    \caption{Structure of a skeletal muscle.}
    \label{fig:3dmuscleskeletalmuscle}
\end{figure}

Thousands of muscle fibres are embedded in a matrix of connective tissue (ECM), which is separated into distinct layers: the epimysium, perimysium, and endomysium (\autoref{fig:3dmuscleskeletalmuscle}). Each of these layers have a complex structure due to the collagen fibres; in particular, the collagen fibres have been observed to have an orientation not necessarily aligned with the muscle fibres \cite{ref:purslow1994}, which gives the ECM an anisotropic response. This gives muscle tissue overall a very complex mechanical response. The fibres render the overall tissue transversely isotropic, while the ECM introduces its own anisotropy. Typically, a simplification can be made to model muscle as overall transversly anisotropic \cite{ref:oomens2009,ref:wakeling2020,ref:zatsiorsky2012} with the anisotropy in the direction of the muscle fibres. The tendon and aponeurosis play an important role in the performance and physiology of muscles, as they connect the muscle to the skeletal system. These tissues are mainly composed of collagen fibres within a matrix of connective tissue. In tendons, the collagen fibres are arranged in parallel, running from the bone to the muscle, while aponeuroses are thin sheets of tendon-like tissue where the collagen fibres are distributed in a more random fashion.

\subsection{One-dimensional Hill-type muscle model}
Most of our current understanding of muscle mechanics comes from measurements in sarcomeres \cite{Moo2020}, single fibres, small muscles (such as rat muscle), and \textit{in silica} experiments of one-dimensional non-linear models. The latter are based on the Hill-type muscle model, introduced by A.V. Hill in \cite{ref:hill1938}. This model provides reasonable estimates of muscle forces during slow contractions \cite{ref:dick2017} with maximal activation of the muscle fibres \cite{ref:perreault2003,ref:wakeling2012,Lu2011,Schmitt2019}. Such models assume that muscles behave as one-dimensional non-linear damped spring systems, and that both mass and inertial effects are negligible; in particular, it assumes that the total force from the muscle is given by the following equation
\begin{align*}
    F_{mus}(\lambda,\epsilon) := F_0(a(t)\hat{F}_{act}(\lambda,\epsilon)+\hat{F}_{pas}(\lambda)),
\end{align*}
where  $\hat{F}_{act}$ is the muscle force due to the activation of the muscle fibres (contractile element force), $\hat{F}_{pas}$ is the non-linear elastic force (passive elastic element), usually referred to as the {\it passive force}, $F_0$ is the maximum isometric force of muscles, $\lambda$ is the one-dimensional stretch, and $\epsilon$ is the rate of change of the stretch over time, commonly referred to as the {\it strain rate}. The function $a = a(t)\in[0,1]$ represents the activation level of the muscle fibres over time. The definitions of the active and passive forces are normalized to the maximum isometric force $F_0$. We comment that in the physiological community and in the context of one-dimensional models, the strain rate $\epsilon$ is generally referred to as the {\it muscle velocity}.

A key assumption in Hill's model is that these forces depend on both the length and velocity at which the fibres are being stretched. The form of the active forces is given by
\begin{align}
    \hat{F}_{act}(\lambda,\epsilon) = \hat{F}_{len}(\lambda)\hat{F}_{vel}(\epsilon),
\end{align}
where $\hat{F}_{len}$ are the forces generated from length changes in the muscle and $\hat{F}_{vel}$ are the forces due to the changes in the strain rate. The force $\hat{F}_{vel}$ satisfies $\hat{F}_{vel}(0) = 1$. The forces $\hat{F}_{len}$ and $\hat{F}_{vel}$ are defined as non-linear functions of the stretch $\lambda$ and the strain rate $\epsilon$ respectively, as seen in \autoref{fig:MaterialResponse}. These are commonly called {\it force-length} and {\it force-velocity} relationships.
Meanwhile, the passive forces $\hat{F}_{pas}$ act to resist fibres from stretching to longer lengths. These play an important role in muscle contraction as non-linear springs: when muscle fibres are lengthened, non-linear passive forces are produced to shorten the fibre back. However, when fibres contract to shorter lengths, experimental data suggest that the passive forces in this regime are negligible, at least for small deformations (\autoref{fig:MaterialResponse}); see, for instance \cite{ref:ross2018-2}.

\subsection{A three-dimensional hyperelastic model for skeletal muscle}
Even though one-dimensional models are able to explain the behavior of whole-muscle mechanics in some specific cases, there are many other cases in which full three-dimensional mechanics is needed to better explain how muscles behave. This is the main focus of the present section.

\subsubsection{Basics of continuum mechanics}
Let us start by recalling some basic concepts in the theory of non-linear elasticity. This overview is mainly based on \cite{ref:holzapfel2001,ref:ogden1997,ref:volokh2016}.
Suppose that an elastic body initially occupies a bounded region $V_0\subset\rrr^3$.
We consider, for fixed $t>0$, that there exists a diffeomorphism $\phi(\cdot,t):V_0\to V$ which maps every point in $V_0$ onto a bounded region
$V\subset\rrr^3$, that is the range of $\phi(\cdot,t)$ is exactly $V$. The region $V_0$ is referred to as the initial (or reference) configuration on the elastic body while $V$ represents the current configuration of such a body. We assume that this map is such that $\phi(\cdot,0)$ is the identity function. Thus, a point $\x\in V$ can be written
as $\x = \phi(\xx,t)$, for some $\xx\in V_0$. We define, at time $t>0$, the {\it material displacement} of a point $\xx\in V_0$ and the {\it spatial displacement} of a point $\x \in V$ as
\begin{align*}
 \uu(\xx,t) := \x - \xx = \phi(\xx,t) - \xx, \quad \u(\x,t) := \x - \xx = x - \phi^{-1}(\x,t).
\end{align*}
For points $\xx \in V_0$ and $\x \in V$, the {\it material velocity} and the {\it spatial velocity} are given respectively by
\begin{align*}
 \vv(\xx,t) := \frac{\partial \phi(\xx,t)}{\partial t} = \Dder{\uu}{t} \equiv \dot{\uu}, \quad \v(\x,t) := \vv(\phi^{-1}(\x,t),t) = \pder{\u}{t}.
\end{align*}
where $\frac{\mathrm{D}}{\mathrm{D}t} \equiv \dot{(\cdot)}$ denotes the \textit{material time derivative} of a material field and $\frac{\partial}{\partial t}$ denotes the {\it spatial time derivative} of a spatial field.
We also define the \textit{deformation gradient tensor} and its determinant (also known as the \textit{Jacobian determinant}) at a given material point to be
\begin{align*}
 \ff(\xx,t) := \ii + \nabla_0\uu(\xx,t), \quad J(\xx,t):= \det\ff(\xx,t), \quad \xx\in V_0,\,\,t>0.
\end{align*}

As for strain measures, we consider the  right and left Cauchy-Green tensors, respectively given by
\[
\cc = \ff^\transpose \ff, \quad \b = \ff \ff^\transpose.
\]
We also introduce the \textit{material velocity gradient} and the \textit{spatial velocity gradient}:
\[
\dot{\ff} = \pder{}{t}\left( \pder{\phi(\xx,t)}{\xx} \right) = \pder{\vv(\xx,t)}{\xx}, \quad \l := \pder{\v(\x,t)}{\xx} = \dot{\ff}\ff^{-1}.
\]
The symmetric part of $\l$ corresponds to what is known as the \textit{rate of strain tensor} $\d$, for which the following holds:
\[
\d := \dfrac{1}{2}\left( \l + \l^\transpose \right), \quad \dot{\cc} = 2\ff^\transpose \d \ff.
\]

Finally, we define the stretch $\lambda$ and strain rate $\epsilon$ of fibres with respect to the reference configuration as
\begin{equation} \label{eq:3dmuscledef-stretch-strain-rate}
\lambda = \|\ff \a_0\|, \quad \epsilon := \dot{\lambda} = \dfrac{(\ff \a_0)^\transpose \, \d \, \ff \a_0}{\lambda},
\end{equation}
where $\a_0 \in \rrr^3$ is a unit vector representing the initial orientation of a fibre at a point $\xx \in V_0$ and $\|\cdot\|$ is the usual 2-norm of vectors in $\rrr^3$.

\subsubsection{Multiplicative decomposition of the deformation tensor} \label{section:multiplicative_decomposition}
In the case of compressible materials (even by small amounts, i.e. nearly incompressible), it is customary to split the deformation tensor into volume-changing (volumetric) and volume-preserving (isochoric) parts, i.e.
\begin{equation}
\ff = J^{1/3} \bar{\ff}, \quad \bar{\ff} := J^{-1/3} \ff.
\end{equation}
We call $\bar{\ff}$ the \textit{modified deformation gradient} and we have that $\det \, \bar\ff = 1$ (this is the isochoric part of the deformation). This decomposition yields modified versions of the right and left Cauchy-Green tensors:
\begin{equation}
\bar{\cc} := \bar{\ff}^\transpose \bar{\ff} = J^{-2/3} \cc, \quad \bar{\b} := \bar{\ff}\bar{\ff}^\transpose = J^{-2/3} \b,
\end{equation}
as well as modified stretch and strain rates:
\begin{equation} \label{eq:modified_stretch_strain_rate}
\bar{\lambda} := J^{-1/3} \lambda, \quad \bar{\epsilon} := \dot{\bar{\lambda}} = \dfrac{(\bar{\ff}\a_0)^\transpose (\ppm : \d) (\bar{\ff}\a_0)}{\bar{\lambda}},
\end{equation}
where $\ppm := \iii - \frac{1}{3} \ii \otimes \ii$ is the spatial projection tensor and $\iii$ is the fourth-order symmetric identity tensor, i.e. in index notation $\iii_{abcd} = \frac{1}{2}(\delta_{ac}\delta_{bd} + \delta_{ad}\delta_{bc})$.

Several invariants (and pseudo-invariants) have been defined for the right Cauchy-Green tensor $\cc$ \cite{ref:weiss1996}. In particular, the models to be presented will make use of the first and fourth modified invariants of $\bar{\cc}$, with derivatives indicated as follows:
\begin{equation}
\bar{I}_1 = \tr \, \bar\cc, \quad \bar{I}_4 = \a_0^\transpose \bar\cc \a_0 = \bar{\lambda}^2, \quad \pder{\bar{I}_1}{\bar\cc} = \ii, \quad \pder{\bar{I}_4}{\bar\cc} = \a_0 \otimes \a_0.
\end{equation}
The fourth invariant $\bar{I}_4$ will be specially important as it contains information on the direction of the fibres. Moreover, its rate of change, $\dot{\bar{I}}_4 =  2\bar\lambda \dot{\bar\lambda}=2\bar{\lambda} \bar{\epsilon}$, depends directly on the strain rate of the fibres.

\subsubsection{A strain-energy function for the muscle-tendon unit}
As described previously, our model of skeletal muscle consists of contributions from the muscle, tendon, fat, base material, and aponeurosis. We consider muscle, aponeurosis, and tendon as anisotropic fibre-reinforced composite and nearly incompressible materials \cite{ref:rahemi2014}. Consequently, the Helmholtz strain-energy function $\Psi$ of the material (as a whole) can be described through a splitting into volumetric and isochoric parts:
\begin{equation} \label{eq:isovolsplit}
\Psi = \Psi_{vol}(J) + \Psi_{iso}(\bar{\cc}, \a_0 \otimes \a_0;\dot{\bar{\cc}}).
\end{equation}
The isochoric part of the strain-energy, $\Psi_{iso}$, encloses the properties of incompressible deformations through the modified right Cauchy-Green tensor (recall that $\det \bar\cc = 1$). Moreover, $\Psi_{iso}$ contains the contributions of all parts of the MTU:
\begin{equation} \label{eq:Psi_iso}
\Psi_{iso} = \Psi_{iso,mus}(\bar{I}_1, \bar{I}_4 ; \dot{\bar{I}}_4) + \Psi_{iso,apo}(\bar{I}_1, \bar{I}_4) + \Psi_{iso,ten}(\bar{I}_1, \bar{I}_4).
\end{equation}
Here, the dependence on $\bar{I}_4$ reflects the fact that the tissues (muscle, aponeurosis, and tendon) are composed mostly of curvilinear fibres running from one point to another within $V_0$. In the current model it is assumed that these fibres are one-dimensional (in the along-fibre direction). However, the space between fibres is not empty and is made of a {\it base material} comprising extracellular material, fluid and other constitutive tissues. The mechanical response of this material include a response from the extracellular material, and the three-dimensional cellular response. This adds a dependence on $\bar{I}_1$ in the isochoric strain-energy above. 
Finally, to include viscoelastic effects in this model, some dependence on $\dot{\cc}$ must be added. We propose to include this rate of deformation as a parameter in the isochoric component of the strain-energy function. 
This can be justified from a numerical analysis perspective: in the resulting system of PDEs, an explicit treatment of the velocity in the time discretization means that the tensor $\dot{\bar\cc}$ (and consequently $\dot{\bar{I}}_4$) is completely known at the current time step. Thus, the stress and elasticity tensors depend on the current displacement (and possibly other current variables) but not on the current velocity. In this way, the fully-dynamic problem can be seen as a sequence of quasi-static deformations, as a pre-computed velocity is used to drive the system from current time step to the next one.

Each one of the contributions in \eqref{eq:Psi_iso} can be further subdivided to consider the different materials found in the MTU, although we save precise definitions for each strain-energy function for \autoref{sec:ParameterFitting}. First, for muscle material, the micro-mechanical effects can be homogenized to include the ECM, cellular material, and fat response. In this work, we utilize the breakdown of the base material derived in \cite{Konno2021} with the addition of a volume fraction encompassing the effects of the fatty tissue:
\begin{equation} \label{eq:homogenizedmuscle}
\begin{split}
\Psi_{iso,mus} &= (1-\beta) \Big\{ \alpha \Psi_{ecm}(\bar{I}_1) + (1-\alpha) \Psi_{cell}(\bar{I}_1)  \\ 
&\hspace{10em} + \Psi_{fibre,mus}(\bar{I}_4;\dot{\bar{I}}_4) \Big\} + \beta \Psi_{fat}(\bar{I}_1),
\end{split}
\end{equation}
where $\alpha$ is the volume fraction of the ECM and $\beta$ is the volume fraction of fat. Then, for the aponeurosis and tendon, we breakdown the ischoric component into base material and fibre components:
\begin{equation}
\Psi_{iso,tis} = \Psi_{fibre,tis}(\bar{I}_4) + \Psi_{base,tis}(\bar{I}_1), \quad tis \in \{apo,ten\}.
\end{equation}
We note here that, while the general formulation of $\Psi_{fibre,tis}$ has a dependency on the strain rate $\bar{\epsilon}$ of the fibres (i.e. on $\dot{\bar{I}}_4$ through \eqref{eq:modified_stretch_strain_rate}), in practice only the muscle fibres will have this dependency. 

\subsubsection{Potential formulation and dynamic equilibrium}

The quasi-static part of the problem can be represented by the stationarity of the following potential function \cite{ref:simo1985}:
\begin{equation} \label{eq:potential}
\Pi(\uu,p,D) = \int_{V_0} \Psi_{vol}(D) + p\left( J(\uu) - D \right) + \Psi_{iso}(\bar\cc) - \int_{V_0} \f_0 \cdot \uu.
\end{equation}
A change in volume during the deformation of the tissues implies a change in internal pressure, which we denote by $p$. The variable $D$ is defined as the true dilation $J(\uu)$ in the reference configuration. Thus, the pressure $p$ plays the role of a Lagrange multiplier.

Notice that \eqref{eq:potential} is defined using quantities in the reference configuration. Thus, a dynamic extension of the Euler-Lagrange equations that arise from the stationarity of this potential, represent a total-Lagrangian description of the system:
\begin{subequations} \label{eq:lagrangian_system}
\begin{gather}
\label{eq:lagrangian_system_1}\dot\uu = \vv, \quad \rho_0 \dot\vv = \bDiv \, \pp(\uu,\vv,p) +  \f_0 \quad \text{in }V_0, \\
\label{eq:lagrangian_system_2}J(\uu)-D = 0, \quad p - \Psi_{vol}'(D) = 0 \quad \text{in }V_0.
\end{gather}
\end{subequations}
where $\bDiv(\cdot)$ is the divergence operator in material coordinates (i.e. $(\bDiv \, \pp)_A = \frac{\partial P_{AB}}{\partial X_B}$), $\rho_0 = \rho_0(\xx)$ is the density of the tissues in the reference configuration, $\f_0$ is a body force, and $\pp$ is the first Piola-Kirchhoff (PK1) tensor. The second Piola-Kirchhoff (PK2) tensor, $\ss$, is related to the PK1 tensor by the identity $\pp = \ff \ss$ and contains the constitutive laws of the material:
\begin{equation}
\ss(\uu,\vv,p) = \ss_{vol}(p) + \ss_{iso}(\uu;\vv) = pJ\cc^{-1} + J^{-2/3} \ppp : \bar{\ss}.
\end{equation}
Here, $\ppp := \iii - \frac{1}{3} \cc^{-1}\otimes \cc$ is the fourth-order material projection tensor and $\bar\ss$ is the ficticious Piola-Kirchhoff stress given by 
\begin{equation} \label{eq:definition_Sbar}
\bar\ss = 2 \pder{\Psi_{iso}(\bar\cc; \dot{\bar\cc})}{\bar\cc} = 2 \left( \pder{\Psi_{iso}(\bar{I}_1)}{\bar{I}_1} \, \pder{\bar{I}_1}{\bar\cc} +  \pder{\Psi_{iso}(\bar{I}_4;\dot{\bar{I}}_4)}{\bar{I}_4} \, \pder{\bar{I}_4}{\bar\cc} \right).
\end{equation}

To reveal the inertial component of the system, we pull \eqref{eq:lagrangian_system} to the current configuration. This constitutes the Eulerian description of the system:
\begin{subequations} \label{eq:3dmuscleeqs-muscle}
\begin{gather}
\label{eq:3dmuscleeqs-muscle1} \pder{\u}{t} = \v, \quad \pder{(\rho \v)}{t} + \bdiv(\rho \v \otimes \v) = \bdiv \, \bsigma(\u,\v,p) + \f \quad \text{in }V, \\
\label{eq:3dmuscleeqs-muscle2}J(\u) - D = 0, \quad p - \Psi_{vol}'(D) = 0 \quad \text{in }V,
\end{gather}
\end{subequations}
where the body force $\f = J^{-1} \f_0$. In particular, inertial forces in the system are represented by the convective term $\bdiv(\rho \v \otimes \v)$. The Cauchy stress tensor $\bsigma$ can be written in terms of the Kirchhoff stress $\btau$ as:
\begin{equation}
J \bsigma =: \btau = \btau_{vol}(p) + \btau_{iso}(\u;\v) = pJ\ii + \ppm : \bar\btau,
\end{equation}
where $\bar\btau := \bar\ff \bar\ss \bar\ff^\transpose$ is the ficticious Kirchhoff stress and $\ppm$ is the spatial projection tensor defined as for \eqref{eq:modified_stretch_strain_rate}.

We close the system given by \eqref{eq:3dmuscleeqs-muscle} (and equivalently \eqref{eq:lagrangian_system}) with initial conditions and mixed boundary conditions. For this, we assume that the boundary in the current configuration $\partial V$ is divided as $\partial V = S_D \cup S_N$, $S_D \cap S_N = \emptyset$. The first part, $S_D$, corresponds to the part of the boundary for which displacements are set (e.g. a clamped face or faces that moves with a prescribed displacement). In turn, $S_N$ represents the set of traction-free surfaces. Mathematically we can write this as
\begin{equation} \label{eq:3dmuscleeqs-muscle3}
\begin{split}
\u = \u_0, \quad \v = \v_0, \quad p = p_0, \quad D = D_0 \quad \text{in }V_0, \\
\u = \u_D \quad \text{on }S_D, \quad \bsigma(\u,\v,p)\n = \mathbf{0} \quad \text{on }S_N.
\end{split}
\end{equation}


\subsubsection{Quasi-static model}
In the quasi-static formulation of the model, we assume that the tissue experiences slow excitations or deformations, on time-scales much slower than that for instantaneous mechanical equilibration. In this case, we do not need the dynamic extension made in \eqref{eq:lagrangian_system}. Therefore, in the Eulerian description, the system \eqref{eq:3dmuscleeqs-muscle} reduces to:
\begin{subequations}\label{eq:3dmuscleeqs-muscle-quasistatic}
\begin{gather}
\label{eq:3dmuscleeqs-muscle1-quasistatic} \bdiv \, \bsigma(\u,p) + \f(\x,t) = \mathbf{0} \quad \text{in }V, \\
\label{eq:3dmuscleeqs-muscle2-quasistatic}D-J(\u) = 0, \quad p - \Psi_{vol}'(D) = 0 \quad \text{in }V, \\
\label{eq:3dmuscleeqs-muscle3-quasistatic}\u = \u_D \quad \text{on }S_D, \quad \bsigma(\u,p)\n = \mathbf{0} \quad \text{on }S_N, \\
\u = \u_0, \quad p = p_0, \quad D=D_0 \quad \text{in }V_0.
\end{gather}
\end{subequations}
These are precisely the Euler-Lagrange equations for the stationarity of the potential \eqref{eq:potential}.

In \autoref{section:3dmusclemethods} we introduce the numerical scheme to approximate the solutions of the non-linear system in \autoref*{eq:3dmuscleeqs-muscle3}. A semi-implicit time stepping is applied at the continuous level to discretize the time variable. The non-linear equations are solved via the application of a Newton's iteration. Finally we discretize the space variables of the resulting linear system with a mixed FEM.

\subsubsection{Muscle mass}\label{section:mass}
We highlight a key question motivating the present work: does muscle mass matter when describing muscle force output and energy? The immediate response may be affirmative. However, these effects may be negligible in certain settings; certainly, in the seminal paper concerning skeletal muscle mechanics \cite{ref:hill1938}, inertial effects due to muscle mass are neglected. These effects are also negligible, to leading order, when dealing with single-fibre experiments, or slow deformations. It is common in the literature to use a quasi-static approximation in these cases. However, it is also known  that Hill-type models (neglecting inertial effect) result in erroneous predictions of actual muscle force in muscles. Recent experimental evidence \cite{ref:holt2014,ross2020} demonstrates the non-negligible impact of muscle mass. In the first of these papers, the authors consider a rat plantaris muscle in which at first only fast-twitch fibres are activated. Next, only slow-twitch fibres are activated. One may expect that when both types of fibres are activated that the maximum velocity of contraction would be {\it slower} than for the fast-twitch types; in fact, it was seen that activating all the fibres lead to an {\it increase} in this velocity. The authors suggest this is because when only one type of fibre is activated, the internal mass of the others leads to a reduction in net work output. In \cite{ross2020}, {\it in-situ}  experiments on rat plantaris muscle were used to study the effects of increasing the mass of muscle on mechanical work. It was seen that the muscles with added mass were able to perform less mechanical work, suggesting again the non-negligible role played by muscle mass in mechanics.

A goal in our work has been to understand if - and under what circumstances - such mass effects need to be included. As discussed earlier, to leading order mass effects are perhaps negligible if we model the entire musculoskeletal system (the rigidity of the skeleton plays a key role), or in single muscle fibres, or at maximal activation. Recent experiments \cite{korta2021} quantify both shear and sound wave speeds in tissue; elastic shear waves are seen to propagate at around 5$ms^{-1}$ (similar to the speed of some neuronal activation pulses), while sound wave speeds are closer to $1.5 \times 10^3 ms^{-1}$. This clearly shows that elastic waves \--- and tissue density \--- cannot be neglected for fast twitches.  

A simple example will illustrate the point. A common experiment in muscle physiology consists of clamping a muscle fibre of initial length $L$ at one end, and pulling it passively with a prescribed displacement $\u_D=g(t)\e_1$ at the other. 
In the quasi-static approximation, all points along the muscle are assumed to respond instantaneously to this displacement, and no compressive elastic wave is expected. In the dynamic model, however, prescribing $\u_D=g(t)\e_1$ becomes equivalent to adding a force proportional to $\frac{X}{L} \ddot g(t)$, where $X$ is the first component of a point in $V_0$. If we, for instance, start to pull the muscle at a constant velocity and then stop doing so, the system experiences an impulsive force. The response will include an elastic wave propagating along the length of the muscle. In \autoref{section:3dmuscleexps}, we contrast the computed mechanical response of the quasi-static and dynamic models.

%% file: Manuscript/discretization.tex
\section{Discretization}\label{section:3dmusclemethods}
In this section we describe the numerical methods used to simulate the three-dimensional deformation of skeletal muscles.  Since the size of the nonlinear system we solve is substantial, we have at this juncture avoided schemes which add to the list of unknowns (by augmenting the system with the unknown velocity) that are solved at each time-step. We note other time-stepping schemes are possible, and are the subject of future work.

\subsection{Eulerian versus Lagrangian descriptions}
At a continuous level, the Lagrangian system \eqref{eq:lagrangian_system} and the Eulerian system \eqref{eq:3dmuscleeqs-muscle} are equivalent, with the interchange between integrals and time derivatives being the link between the two descriptions \cite{ref:holzapfel2001}. Therefore, when the time variable is discretized first, these two descriptions are no longer equivalent, but because they are descriptions of the same physical process, the end result is very similar. 

To avoid the complications brought by the convective terms in the Eulerian description, we choose to discretize the Lagrangian description \eqref{eq:lagrangian_system}. Notice also that, unlike in the Eulerian description where the domain $V$ changes at every time step, the Lagrangian discretization only requires a mesh of the reference configuration $V_0$ that is kept fixed throughout the entire simulation.

\subsection{Total Lagrangian formulation}
To approximate a solution of the non-linear dynamic system in \eqref{eq:3dmuscleeqs-muscle}, we use a semi-implicit time stepping for the time derivatives, and a conforming finite element scheme for the spatial variables. The displacement,
dilation, and pressure are updated with an implicit approach while the velocity is computed with an explicit update. The velocity is treated explicitly and it is only updated once the new displacement is computed. Using this time stepping in \eqref{eq:lagrangian_system}, we propose the problem: at time step $t_n$, find a material displacement $\u^n$ (we use lowercase letters from now on for a more pleasant reading), a material pressure $p^n$, and a material dilation $D^n$ such that
\[
\left\{\begin{gathered}
\rho_0 \u^n - \delta t^2 \bDiv \, \pp(\u^n, \v^{n-1}, p^n) = \rho_0 \u^{n-1} + \delta t \rho_0 \v^{n-1} + \delta t^2 \f_0^n, \\
J(\u^n) - D^n = 0, \quad p^n - \Psi_{vol}'(D^n) = 0,
\end{gathered} \right.
\]
where $\delta t > 0$ is the time step size. The first equation comes from combining the two equations in \eqref{eq:lagrangian_system_1} using the time discretization $\v^n = \left( \u^n - \u^{n-1} \right)/\delta t$, which is also used to compute the material velocity as a post-process at each time step.

Multiplying these equations by a test function $\delta \bxi := (\delta\u,\delta p, \delta D) \in \mathcal{X}$, with $\mathcal{X} := \hh^1(V_0) \times L^2(V_0) \times L^2(V_0)$, and integrating by parts one obtains the following non-linear variational formulation: given the previous displacement $\u^{n-1}$, previous velocity $\v^{n-1}$, force $\f_0^n$, and a prescribed displacement $\u_D^n$, we are to find $\bxi^n := (\u^n, p^n, D^n) \in \mathcal{X}$ such that $\u^n = \u^n_D$ on $S_D$ and
\begin{equation} \label{eq:def_R}
R(\bxi^n , \delta\bxi) := R_{ine}(\u^n,\delta\u) + R_{int}(\bxi^n,\delta\bxi) - R_{ext}(\u,\delta\u) = 0, \quad \forall \, \delta\bxi \in \mathcal{X},
\end{equation}
that is, for all virtual displacements $\delta \u$, virtual pressures $\delta p$, and virtual dilations $\delta D$. Such a formulation can be understood as a semi-discrete version of the principle of virtual work in the total Lagrangian formulation. The semi-discrete work due to inertial, internal, and external forces are given respectively by:
\begin{subequations}
\begin{align}
\begin{split}
&R_{ine}(\u^n,\delta \u) =  \delta t^{-2}  \iprod{\rho_0 \u^n , \delta\u} \\
&\hspace{8em}- \delta t^{-2} \iprod{\rho_0 \u^{n-1}, \delta \u} - \delta t^{-1} \iprod{\rho_0 \v^{n-1} , \delta \u},
\end{split} \\
\begin{split}
&R_{int}(\bxi^n, \delta \bxi) = \iprod{\btau(\u^n,\v^{n-1},p^n) , \nabla_0(\delta\u) \ff(\u^n)^{-1}}  \\
&\hspace{8em}+ \iprod{J(\u^n) - D^n, \delta p} + \iprod{\Psi_{vol}'(D^n)-p^n, \delta D},
\end{split} \\
\begin{split}
&R_{ext}(\u^n,\delta\u) = \iprod{\f_0^n, \delta \u}.
\end{split}
\end{align}
\end{subequations}
The non-linear nature of this problem is handled using Newton's method. At each time step, given $\u^{n-1}, \v^{n-1}, \u_k^n, p_k^n, D_k^n$, we want to find an update $\d\bxi_k := (\d\u_k, dp_k, dD_k)$ such that $\d\u_0 = \u^n - \u^{n-1}$ on $S_D$, $\d\u_k = \zero$ on $S_D$ for $k \geq 1$, and 
\begin{equation}
\left[ \mathcal{D} R (\bxi^n_k)(\d\bxi_k), \delta\bxi \right] = -\left[ R(\bxi_k^n),\delta \bxi \right],
\end{equation}
for all virtual increments $\delta \bxi \in \mathcal{X}$, to then update the solution as $\bxi_{k+1}^n = \bxi_k^n + \d\bxi_k$, for all $k \geq 0$ until convergence. At the current time step, we set the initial solution for the Newton iteration as the converged solution from the previous step, i.e. $\bxi_0^n = \bxi^{n-1}$. The tangent operator, $\mathcal{D}R$, corresponds to the Gateaux derivative of the operator $R$ defined in \eqref{eq:def_R}, and $[\cdot,\cdot]$ denotes the duality pairing between $\mathcal{X}$ and $\mathcal{X}^*$. For more information on the linearization of this operator and the resulting fourth-order elasticity tensors, we refer the reader to \cite{ref:holzapfel2001,ref:wriggers2009}.

\subsection{Finite element discretization}
To discretize the linear operator $\mathcal{D}R(\bxi_k^n)$, we use a standard conforming finite element. Let $\tttt_h$ be a regular mesh of $\overline{V_0}$ made of (isoparametric) brick elements $T$ of diameter $h_T$. The meshsize is simply defined as $h := \max \left\{ h_T : T \in \tttt_h \right\}$. We use the same finite element spaces as those used in \cite{ref:pelteret2012}, that is, $\qq_{k+1} \times P_k \times P_k$. For a non-negative integer $k$ and an element $T$ of $\tttt_h$, let us consider the $P_k(T)$ as the space of all polynomials of degree at most $k$ defined in $T$, with $\pp_k(T)$ denoting the vector version of $P_k(T)$. Then
\[
\begin{aligned}
\qq_{k+1} &:= \left\{ \w \in \mathcal{C}(\overline{V_0}) \, : \, \w\big|_T \in \pp_{k+1}(T), \ \forall \, T \in \tttt_h \right\}, \\
P_k &:= \left\{ q \in L^2(V_0) \, : \, q\big|_T \in P_k(T), \ \forall \, T \in \tttt_h \right\}.
\end{aligned}
\]
Normally, we will use $k=1$, that is, displacements are approximated using quadratic 27-node Lagrange elements $\qq_2$, while pressures and dilations are approximated using linear discontinuous elements $P_1$ (based on monomials, 4 degrees of freedom per element). We also note that the element $\qq_1 \times P_0 \times P_0$ does not exhibit locking though the theory is lacking.

%% file: Manuscript/strainenergies.tex
\section{Determination of model parameters and strain-energy description} \label{sec:ParameterFitting}
Here we describe the methods used to determine the parameters in the material response of the muscle tissues. The full list of parameters and functions are included in the \textcolor{red}{Supplementary Material}.

\subsection{Volumetric response}\label{section:volumetricresponse}
In this study we adopt the form of the volumetric stress used in \cite{ref:pelteret2012,ref:rahemi2014} (see also \cite{ref:holzapfel2001} for more details on volume changing free energy functions), which is typically used for nearly incompressible biological materials \cite{ref:holzapfel2001}:
\[
\Psi_{vol}(J) = \dfrac{\kappa}{4} \left( J^2 - 2 \ln \, J - 1 \right), \quad \kappa(\xx) = \sum_{tis \in \{mus,apo,ten\}} \kappa_{tis} \bchi_{tis}(\xx),
\]
in which $\kappa_{tis}>0$ is the bulk modulus of the corresponding tissue and $\bchi_{tis}$ is the indicator function of the tissue subdomain. In particular, the muscle bulk modulus is computed as $\kappa_{mus} = (1-\beta) \left(\alpha \kappa_{ecm} + (1-\alpha) \kappa_{cell}\right) + \beta \kappa_{fat}$. The values for bulk moduli are determined through ensuring that the near incompressibility condition is met, ie. the values for the moduli are chosen so that there is little change in the volume of the muscle. Based on \cite{ref:blemker2005,rahemi2014regionalizing}, we choose $\kappa_{ecm} = 1\times 10^6$Pa. Due to the incompressibility of water and the fact cells are mainly composed of water we choose the volume fraction for the cellular material and fat to be $\kappa_{cell} = \kappa_{fat} = 1\times10^7 $Pa \cite{Konno2021}. It has been shown in \cite{ref:gardiner2001,Konno2021} that the resulting strain mechanics are not very sensitive to the bulk moduli, so we mainly choose these values to ensure that the volume change is small. Due to the increased stiffness of the aponeurosis and tendon materials, the bulk modulus is chosen to be larger, $\kappa_{apo} = \kappa_{ten} = 1\times 10^8$Pa.

\subsection{Isochoric response due to base material properties}
The mechanical responses of the isochoric components of the base material for each of the different tissues  (including ECM and cellular components for muscle) are modeled using Yeoh-type strain-energy functions \cite{ref:yeoh1993}. In turn, the intramuscular fat was modeled as a material of Neo-Hookean type. This yields two types of strain-energies that will be used in this work:
\begin{equation}
\Psi_{Yeoh}(\bar{I}_1) = \sum_{k=0}^3 c_k \left( \bar{I}_1 - 3 \right)^k, \quad \Psi_{Neo}(\bar{I}_1) = c_1 \left( \bar{I}_1 - 3 \right).
\end{equation}
where the constants $c_k \in \rrr$, $k=0,\dots,3$, are determined based on experimental data for each one of the tissues. For the Yeoh model, the basis of the fitting procedure will consider the incompressible version of the Cauchy stress response:
\begin{equation} \label{eq:constitutiverel}
\bsigma_{Yeoh} = -p_{Yeoh} \ii + 2 \pder{\Psi_{Yeoh}(I_1)}{I_1} \b, \quad p_{Yeoh} := \dfrac{2}{\lambda} \pder{\Psi_{Yeoh}(I_1)}{I_1},
\end{equation}
where $I_1 := \tr \, \b$ is the first invariant of the left Cauchy-Green tensor.

In our model, any deformation transverse to the fibre direction will be largely due to the base material. Hence we determine the parameters, $c_k$, using experimental data for tissue pulled in the direction transverse to the fibres. Experimental data is only available for the stress as a function of stretch in a material; therefore, we use the constitutive relationship \eqref{eq:constitutiverel} and express the left Cauchy-Green tensor $\b$ in principal stretches. In particular, we assume that for uniaxial stretch the principal components of the stretch, $\lambda_i\in\rrr$, $i=1,\ldots,3$, are given by $\lambda = \lambda_1$ and $\lambda_2 = \lambda_3$, here $\lambda_1$ is the direction of the stretch applied in the experiment.
These relations along with $\lambda_1 \lambda_2 \lambda_3 = 1$ (which comes from the incompressibility assumption) gives us the stress response from the Yeoh model in the principal spatial direction as
\begin{align}
    \sigma_{Yeoh}(\lambda) = 2 \left(\lambda^2 - \frac{1}{\lambda}\right) \left[c_1 + 2c_2 \left(\lambda^2 + \frac{2}{\lambda}- 3\right) + 3c_3 \left(\lambda^2 + \frac{2}{\lambda} - 3\right)\right].
\end{align}
The parameters, $c_k$, were then obtained via non-linear regression to data from \cite{ref:mohammadkhah2016,Gillies2011decell,Alkhouli2013,Jin2013} to determine the muscle material constants and \cite{ref:mechazizi2009} for the aponeurosis and tendon base material. For the ECM component we utilize data from a decellularized ECM, which, after decellularization, undergoes material testing \cite{Gillies2011decell}. For the cellular component of the base material, obtaining mechanical data is more difficult \cite{Konno2021}. The cellular component consists largely of different cells, e.g. muscle fibre cells, satellite cells, and neuron cell bodies; thus, we can utilize material data obtained from brain grey matter, which is a mass of cells with little connective tissue \cite{Jin2013}. The fat data is obtained through tensile experiments on subcutaneous fat \cite{Alkhouli2013}. This type of fat is chosen, as it has little anisotropy in the material response, which is ideal as we model the base material and therefore the fat response as isotropic. Additionally, this was the only fat data to our knowledge that contained a tensile response and was in a region near skeletal muscle, as many other studies investigated the compressive response of fat.
The parameter values from the regression are given in the \textcolor{red}{Supplementary material}.

\subsection{Isochoric response due to along-fibre properties}
\begin{figure}[!ht]
    \centering
    \includegraphics[width = 0.3\textwidth]{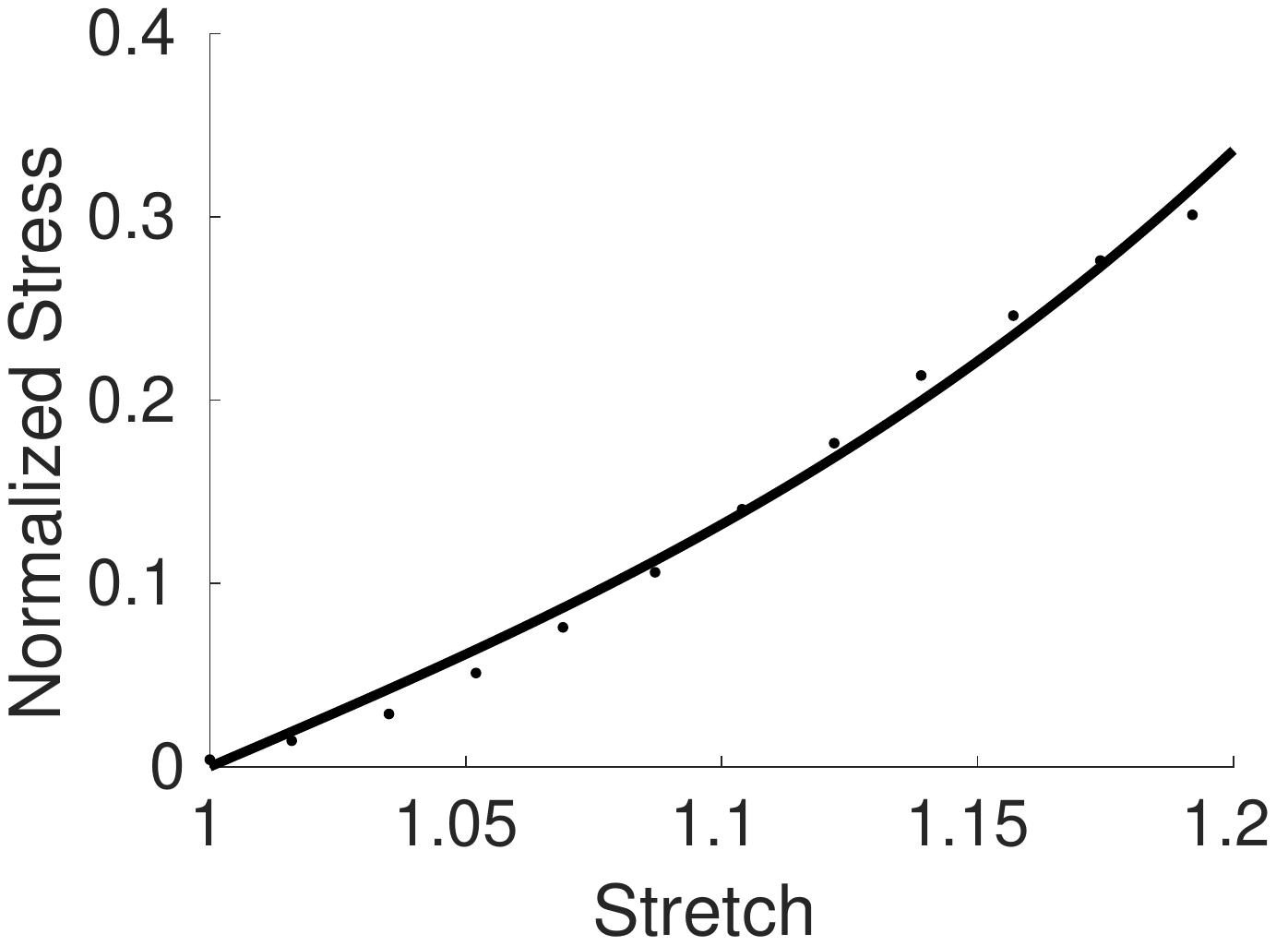}
    \includegraphics[width = 0.3\textwidth]{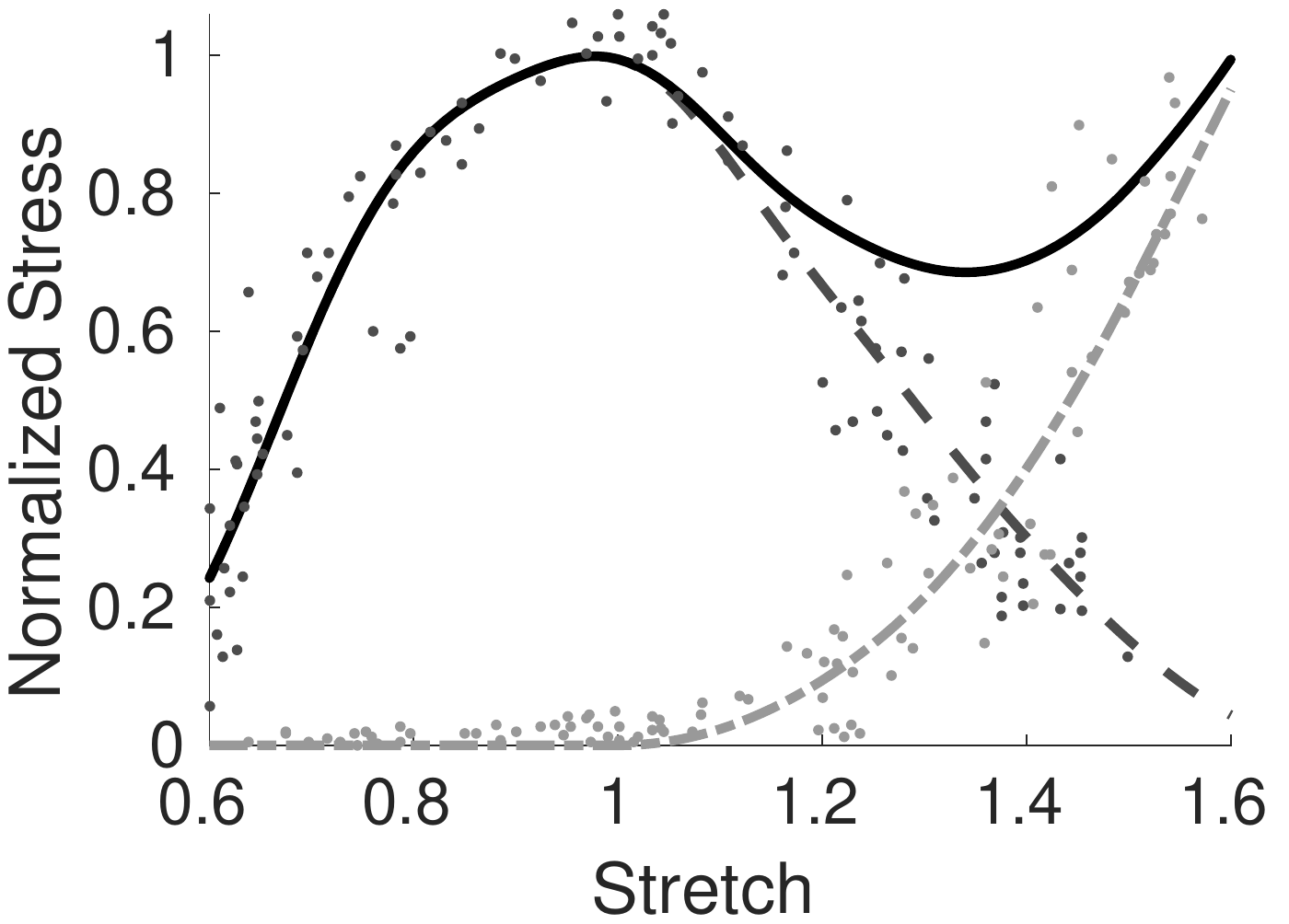}
    \includegraphics[width = 0.3\textwidth]{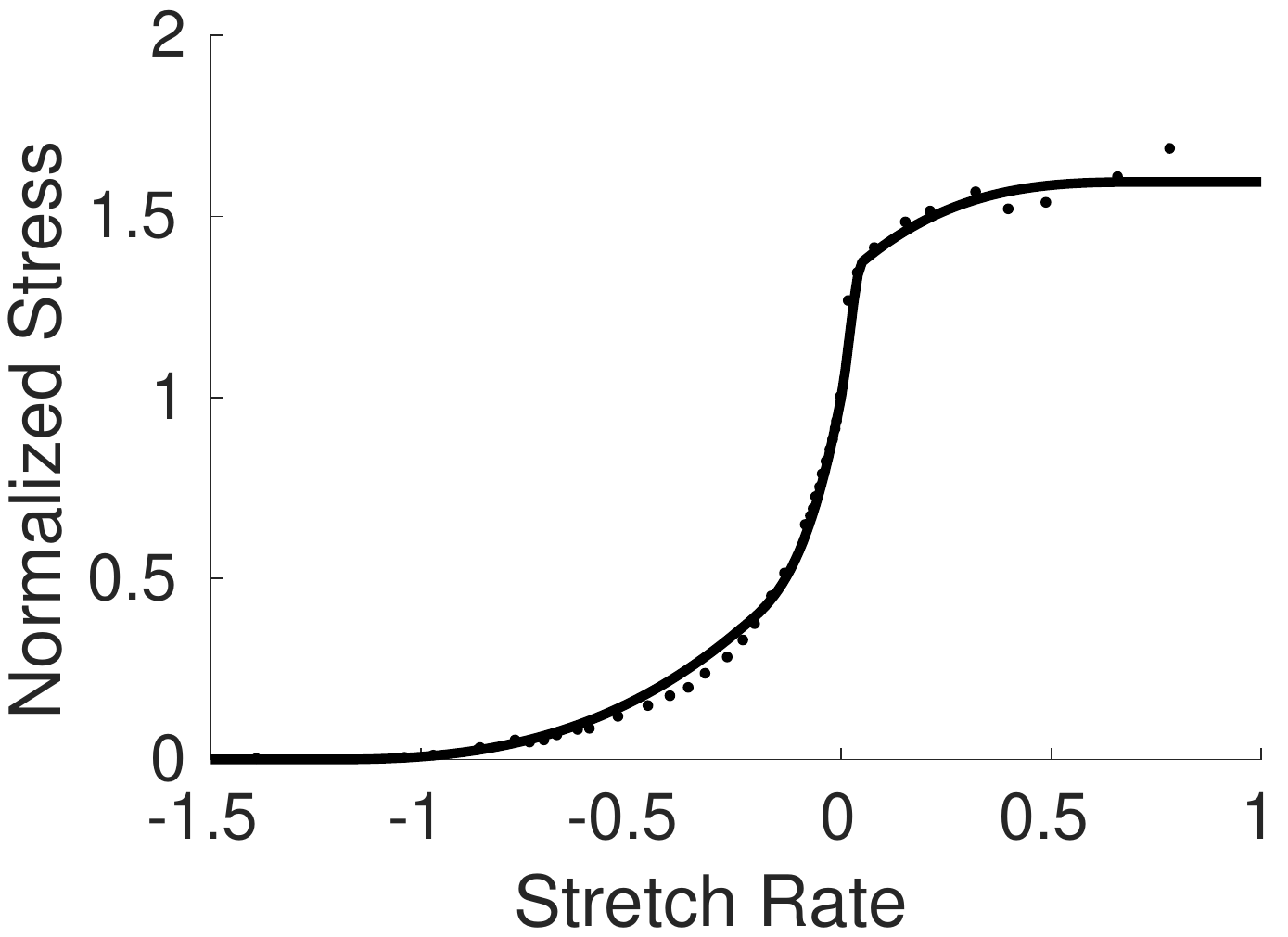}
    \caption{Skeletal muscle material response for the isochoric component of the base material fit to data from \cite{ref:mohammadkhah2016} (A). The force-length active and passive response is fit to data from \cite{ref:winters2011} (B). The response for the force-velocity component is fit to data from \cite{Roots2007} (C).}
    \label{fig:MaterialResponse}
\end{figure}

For the fibre components of the strain-energy function (i.e. $\Psi_{fibre,tis}$, $tis \in \{mus,apo,ten\}$), we use an implicit model where the constitutive laws are defined in terms of derivatives of these strain-energies (which is actually what is needed to construct the fictitious Piola-Kirchhoff stress, $\bar\ss$, see \eqref{eq:definition_Sbar}):
\begin{equation}
2\bar{I}_4 \pder{\Psi_{fibre,mus}(\bar{I}_4;\dot{\bar{I}}_4)}{\bar{I}_4} = \sigma_{mus}(\bar\lambda,\bar\epsilon), \quad 2\bar{I}_4 \pder{\Psi_{fibre,tis}(\bar{I}_4)}{\bar{I}_4} = \sigma_{tis}(\bar\lambda),
\end{equation}
for $tis \in \{apo,ten\}$,  with $\bar{I}_4 = \bar{\lambda}^2$ and $\dot{\bar{I}}_4 = 2\bar{\lambda} \bar{\epsilon}$ (see \autoref{section:multiplicative_decomposition}). Similar to the formulation of Hill's model, we can decompose the response from $\sigma_{mus}$ as
\begin{equation}
\sigma_{mus}(\bar{\lambda}, \bar{\epsilon}) = \sigma_0 \left( a(\xx,t) \, \hat{\sigma}_{act}(\bar{\lambda},\bar{\epsilon}) + \hat{\sigma}_{pas}(\bar{\lambda}) \right),
\end{equation}
where $\sigma_0 > 0$ is the maximum isometric stress in the muscle fibres, $a(\xx,t) \in [0,1]$ is the fibre activation profile, $\hat{\sigma}_{act}$ is the active stress, and $\hat{\sigma}_{pas}$  is the passive stress in the fibre. In general, the active stress is defined by two forces of different nature: elastic forces due to the stretch of the muscle fibres and forces due to the velocity at which the muscle fibres are being stretched. The active stress is defined as
\[
\hat{\sigma}_{act}(\bar{\lambda}, \bar{\epsilon}) :=  \hat{\sigma}_{len}(\bar{\lambda})  \, \hat{\sigma}_{vel}(\bar\epsilon).
\]
Here, $\hat\sigma_{len}$ is called the {\it force-length} function and $\hat\sigma_{vel}$ is the {\it force-velocity} function (see \autoref{fig:MaterialResponse}). For the sake of clarity, in this case the ``velocity'' refers to the velocity in which the fibres are stretched, i.e. the strain rate $\bar\epsilon$ (which should not be confused with either the material or spatial velocity, see \eqref{eq:modified_stretch_strain_rate}). We also note that the stresses $\hat\sigma_{len}$ and $\hat\sigma_{vel}$ are normalized so that $\hat\sigma_{len}(1) = 1$ and $\hat\sigma_{vel}(0) = 1$. In particular, for the quasi-static setting in \eqref{eq:3dmuscleeqs-muscle-quasistatic}, we have $\v=\zero$, and so $\hat{\sigma}_{act} = \hat{\sigma}_{act}(\bar{\lambda}) = a(\xx,t) \, \hat{\sigma}_{len}(\bar{\lambda})$.

In contrast, aponeurosis and tendon are not subject to active forces, so the stresses take a simpler form:
\[
\sigma_{tis}(\bar{\lambda}) = \sigma_{0,tis} \hat\sigma_{len,tis}(\bar\lambda), \quad tis \in \{apo,ten\},
\]
where $\sigma_{0,tis} > 0$ are the maximum isometric stresses and $\hat\sigma_{len,tis}$ are force-length functions of each tissue.

\subsection{Parameter fitting for force-length and force-velocity functions}
Many different functions have been used to represent the force-length and velocity relationships \cite{Rockenfeller2017}, including trigonometric (e.g. \cite{Scott1991}) and polynomial functions (e.g. \cite{He1991}); however, not all these functions are equal, some will be more computationally efficient than others \cite{Rockenfeller2017}. A balanced option developed in previous work \cite{ref:ross2018-1} uses cubic B\'ezier curves to capture the intrinsic properties of muscle. B\'ezier curves are determined based on a set of polynomials that smoothly interpolate user defined points \cite{PrautzschHartmut2002BaBt}. In our case the user defined points are chosen to achieve a smooth function that accurately captures the intrinsic muscle properties. A cubic polynomial is used over higher order polynomials to allow for better local control of the function.

These B\'ezier curves have been used in previous one-dimensional modelling experiments to investigate cyclic work loops \cite{ref:ross2018-1}; however, in our three-dimensional model of muscle, they are not as practical. We opted to use a along fibre response that was twice continuously differentiable. The reason for this is that that in practice more computational time is required to compute the resulting stress with the B\'ezier curve. A cubic polynomial was fit to the existing B\'ezier curves from \cite{ref:ross2018-1} for the passive force-length and force-velocity response, but with the additional requirement that the function is twice continuously differentiable everywhere. For the active force-length relationship a trigonometric function was used opposed to the polynomial fit, but with the same constraint that the function is twice continuously differentiable everywhere.

%% file: Manuscript/experiments.tex
\section{Numerical experiments}\label{section:3dmuscleexps}
We now present a number computational studies based on our model, exploring the impact of muscle mass on dynamics, as well as probing the role of structural changes during cerebral palsy. In \autoref{sec:DynamicMuscle} we demonstrate the dynamical behavior of the model, contrasting with the use of a quasi-static approximation. Furthermore, in \autoref{sec:CPmuscleExperiments} we show an application of the quasi-static model to muscle effected by cerebral palsy, which demonstrates the ability of the model to capture the influence of muscular disorders. Our framework allows us to include patient-specific meshes based on MRI data, and an experiment demonstrating this is reported in \autoref{sec:mir-geometry}.  All experiments have been implemented in the finite element library deal.ii version 8.5.1 \cite{ref:davydov2017}. In each, we use mapped hexahedral elements with smallest $\Delta x= .75$mm, and  a time step of $10^{-5}mm$. 

\subsection{Initial and boundary conditions}\label{section:3dmuscleinitial}
In what follows, we assume that at time $t=0$, the muscle is at rest, unactivated, and in an underformed state. Therefore $\u(\xx,0)=\mathbf{0}$, $\v(\xx,0) = 0$ and $\ff(\xx,0) = \ii$. The latter implies $D(\xx,0) = 1$ and
\begin{align*}
    p(\xx,0) = \Psi_{vol}'(D)\big|_{t=0} = \frac{\kappa}{2}\left.\left(D-\frac{1}{D}\right)\right|_{t=0} = 0.
\end{align*}
Note that this choice of initial conditions is compatible with the definitions introduced in \autoref{section:3dmusclemechanics}.

\begin{figure}[!ht]
    \centering
    \includegraphics[width=0.69\textwidth]{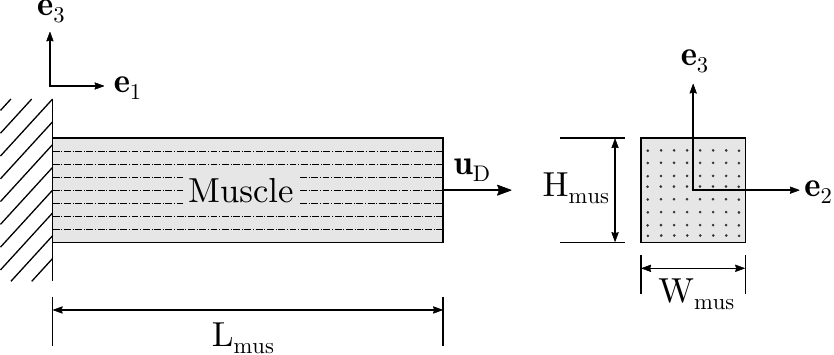}\hspace{1em}\raisebox{1.68cm}{\small \begin{tabular}{|c|c|}\hline
    $L_{mus}$ & 52.0008 [mm] \\\hline
    $W_{mus}$ & 13.7500 [mm] \\\hline
    $H_{mus}$ & 5.5783 [mm] \\\hline
\end{tabular}}\\
    {\includegraphics[width=0.69\textwidth]{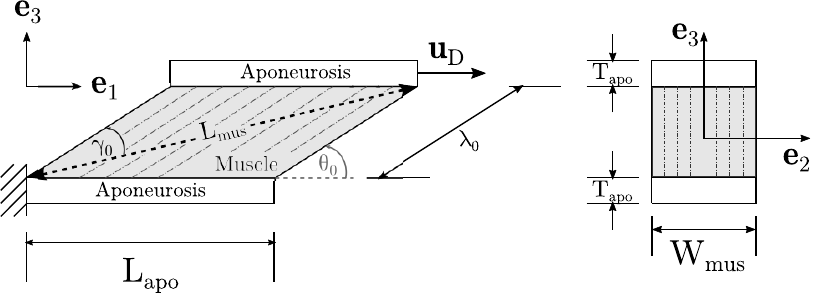}}\hspace{1em}\raisebox{1.5cm}{\small \begin{tabular}{|c|c|}\hline
    $L_{apo}$ & 52.0008 [mm] \\\hline
    $\lambda_0$ & 16.2500 [mm] \\\hline
    $\theta_0$ & 20$^\circ$ \\\hline
    $L_{mus}$ & 67.5000 [mm] \\\hline
    $T_{apo}$ & 0.7500 [mm] \\\hline
    $W_{mus}$ & 13.7500 [mm] \\\hline
    \end{tabular}}
    \caption{Geometries considered for the numerical experiments: a block of pure muscle (top) and the simplified version of the medium gastrocnemius from \cite{ref:Ross2021} (bottom). Note that in the block geometry, the fibres run parallel to the $x$-axis, whereas in the medium gastrocnemius they run at an angle $\theta_0$ from the $x$-axis. In this case, the relationship between the length of the muscle $L_{mus}$ and the initial fibre length $\lambda_0$ is given by $\lambda_0 \sin(\theta_0) = L_{mus}\sin(\theta_0 - \gamma_0)$.\label{fig:geometries}}
\end{figure}

The boundary conditions that will be implemented in the following numerical experiments are the same with exception of the exact form of the Dirichlet condition applied to the muscle. \autoref{fig:geometries} shows diagrams of the geometries in consideration. We prescribe a zero Dirichlet condition on the $-x$-face of the geometries (i.e. $\u_D = \zero$). In turn, the condition for the $+x$-face has the form $\u_D(\cdot,t) = \left(d(t), 0,0\right)^\transpose$, with $d(t)$ being a smooth function that varies according to the experiment. The rest of the boundary correspond to traction-free surfaces.

\subsection{Dynamic muscle mechanics}\label{sec:DynamicMuscle}
We begin by discussing the simple experiment discussed in \autoref{section:mass}. Consider a block of pure muscle, clamped at one end with fibres oriented along the x-axis, as in \autoref{fig:geometries}. The length of the muscle block is initially $L_{\mathrm mus}$, and points on the $x_{end}$ face of of block are subject to a Dirichlet boundary condition $\u_d= (d(t),0,0)$, where 
$$ 
d(t)=\begin{cases}  0,& 0\leq t \leq 0.05\\ (t-0.05)\cdot0.1L_{\mathrm mus}, & 0.05 \leq t \leq .015   \\0.1L_{\mathrm mus}& t\geq 0.15.\end{cases} 
$$
On the transverse faces, a zero traction is prescribed.  We simulate the elastic response in the muscle block over time using a $\qq_2$-$P_1$-$P_1$ discretization on a fine mesh. 

\begin{figure}[!ht]
    \includegraphics[width=\textwidth]{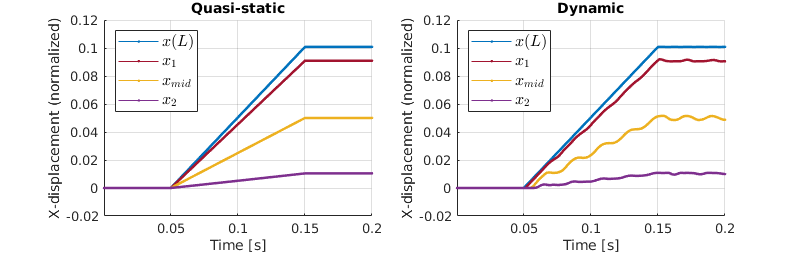}\vspace{-1em}
    \caption{The x-component of displacement as a function of time, at different locations along a parallel-fibered muscle block. One face is fixed, and the opposite one is initially held at a fixed distance, then pulled with a constant velocity, and then held fixed again. Right: Quasi-static calculations. Left: Fully dynamic calculations. Elastic waves are seen clearly propagating in the dynamic model. The normalization factor is $L_{mus}$ (see \autoref{fig:geometries}).}\label{fig:dyn-isocontrast}
\end{figure}

In \autoref{fig:dyn-isocontrast}, we present time traces of the x-components of the displacement. These are measured at 4 distinct points along the centerline of the muscle block: at $x(L)$ (the face being pulled), $x_1$ close to this face, $x_{mid}$ at the midpoint of the block and $x_2$ close to the clamped face. The left plot in \autoref{fig:dyn-isocontrast} shows these traces when we use a quasistatic approximation, and on the right when we use the full dynamic model. We note that, as expected, a quasistatic approximation leads to all points along the muscle reacting instantaneously to the pulling (albeit with differing magnitudes). On the other hand, inertial effects are clearly visible in the dynamic calculation; an elastic wave propagates towards the clamped end. 

\subsection{Isokinetic force-velocity experiments}
The force-velocity relationship is a fundamental principle of skeletal muscle physiology based on Hill's studies on isolated frog muscles (\cite{Alcazaretal}). In this study, Hill found that when a muscle is stimulated isometrically and then suddenly released under a load, the shortening speed increases with decreasing loads in a hyperbolic relationship (\cite{ref:hill1938}). While Hill's experiments have been instrumental to our understanding of the muscle's force-velocity relationship, due to the small scale of these experiments, important factors such as mass and the three-dimensional behavior of muscle tissue were neglected. For instance, consider  single-fibre experiments (or those on small muscles) where one end is clamped, and the other is moving at a fixed applied velocity of $\v_{appl}$. One reasonable assumption may be that each point along a fibre experiences a stretching at a rate given by $\v_{appl}\cdot \a$, where $\a=\bar\ff\a_0$ is the local direction of the fibre. We term this an {\it isokinetic} assumption; here, mass is included in the dynamic model. Another reasonable assumption may be that it suffices to consider a quasi-static approximation (without mass), and in this case the free end moves with $\u_D=\delta t\,\v_{appl}$.
In a fully dynamic calculation, the local fibre strain rate is given by \eqref{eq:3dmuscledef-stretch-strain-rate}. 

In this experiment, we compare these three models. All computations are performed on a sample medium gastrocnemius (bottom figure in  \autoref{fig:geometries}), in which the fibres are pennate (i.e., not parallel to the direction of the aponeurosis). The configuration is held fixed along the $-x$-face of the bottom aponeurosis. The $+x$-face of the top aponeurosis is initially held fixed while the muscle fibres are activated. Then this face is subject to a deformation $\u_D = -\v_{appl}\cdot\e_1$, corresponding to a shortening.

\begin{figure}
	 \centering
	\includegraphics[width=0.9\textwidth]{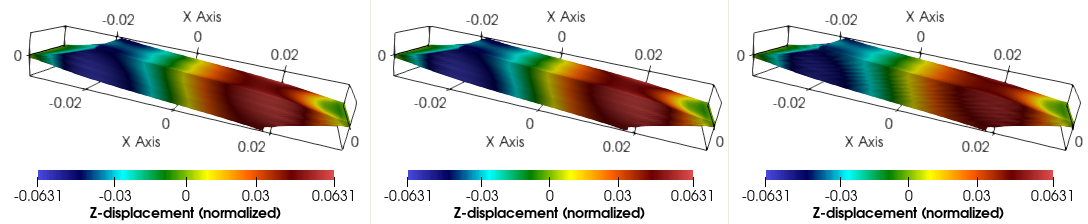}\vspace{0.5em}\\
	\includegraphics[width=0.9\textwidth]{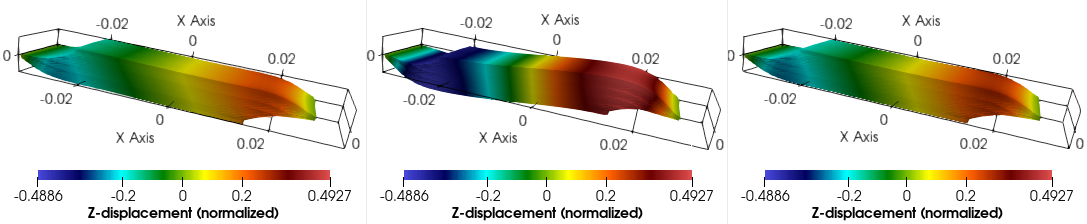}
 	\caption{Left to right: A quasi-static, iso-kinetic and fully dynamic simulation on a model gastrocnemius. Top row: 45\% activation. Bottom row: Activated muscle after shortening phase. The normalization factor corresponds to the total height of the unit, i.e. $\lambda_0 \sin(\theta_0) + 2 T_{apo}$ (see \autoref{fig:geometries}).}\label{fig:threemodels}
 \end{figure}

In each of these simulations, the muscle is activated initially, and then shortened. In \autoref{fig:threemodels}, we show the deformed muscles; the colors indicate the (normalized) z-component of displacement. All three computations yield similar results during the initial activation phase (i.e., for a low activation level), but substantially different behaviors once the muscle is shortening (bottom 3 figures in \autoref{fig:threemodels}). Once mass effects are included (in the isokinetic and the dynamic models), the muscle is seen to curve more towards the end where the shortening is being applied; both the magnitude of this curvature as well as the pattern of deformation are clearly not the same.

\subsection{Quasi-static muscle mechanics: modelling cerebral palsy}\label{sec:CPmuscleExperiments}
The mechanics of skeletal muscle can be substantially altered by muscular disorders, aging, and disuse \cite{ref:lieber2013,Howard2021}. The altered mechanics come as a result of changes to the microstructure of the muscle, such as increased extracellular-matrix (ECM) volume fraction, and changes to the overall muscle architecture. Understanding how individual changes affect muscle function is difficult to do using experimental techniques; however, utilizing a modelling approach we can understand the direct consequence of changes occurring during muscular disorders. In this work we demonstrate an example of cerebral palsy (CP) affected muscle.  To consider the effects from the micro-mechanical effects of CP, we vary the volume fraction and stiffness parameters in \eqref{eq:homogenizedmuscle}. Another characteristic of CP affected muscle is increased length of the sarcomeres; mathematically, this can be modelled as a shift in the fibre force-length relationship $\sigma_{mus}(\bar\lambda + c_{\text{sarco}},\bar\epsilon)$ \cite{Konno2022}.
 
To evaluate the influence of the CP properties in skeletal muscle, we can perform a force-length test, using the quasi-static model of muscle (see e.g. \cite{Lieber2004}), aiming to mimic typical force-length experiments holding the end of a muscle fixed. The goal is to study the impact of various structural changes to a passively-stretched muscle. To construct a passive force length curve, the muscle is pulled to a given length. The magnitude of the force is then measured on the end $+x$-face of the muscle. The active force-length curve can similarly be created by pulling the muscle to a given length, then activating using a linear activation time ramp to 100\%.

\begin{figure}[!ht]
    \centering
    \includegraphics[width=0.95\textwidth]{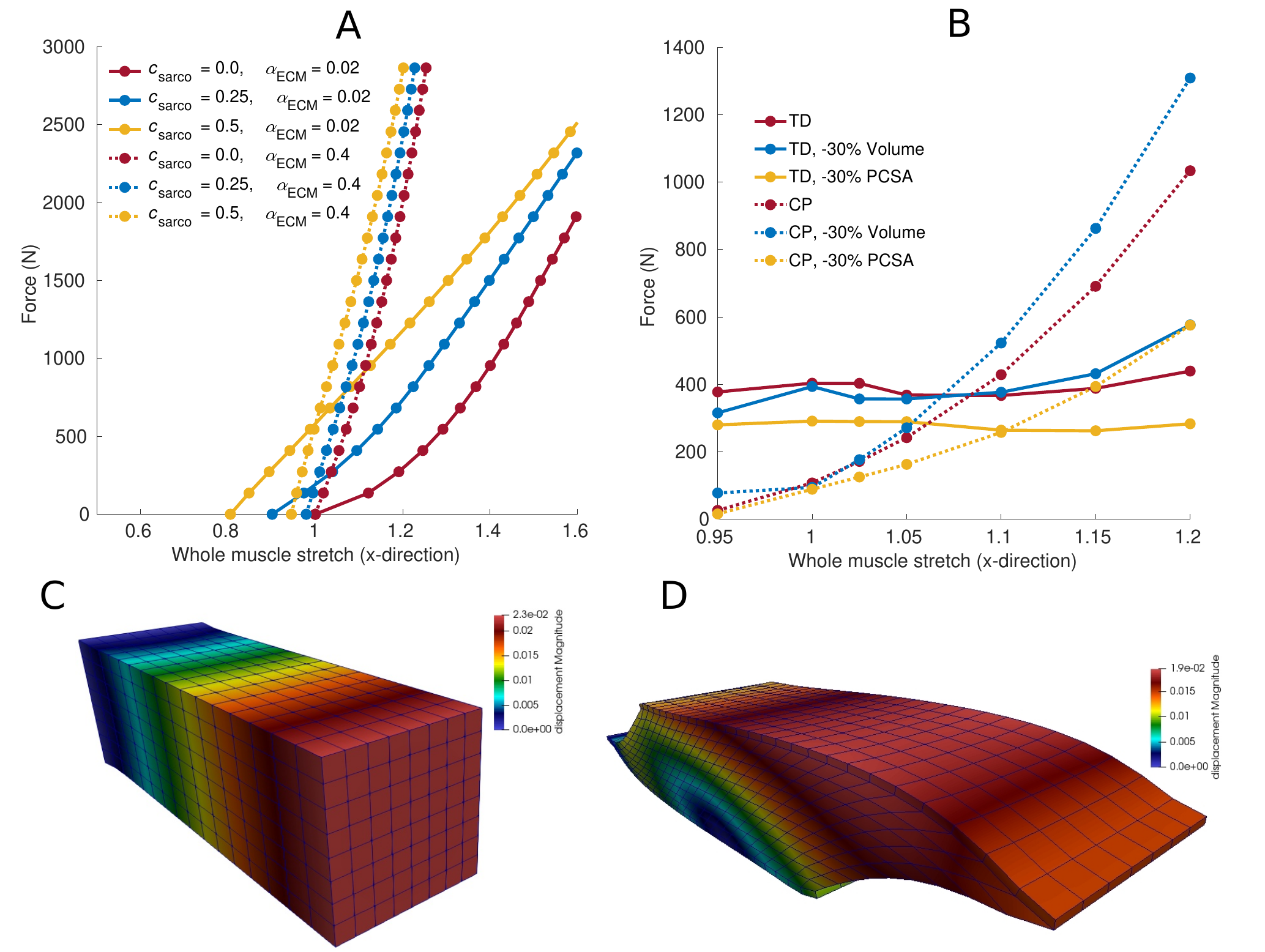}
    \caption{Modelling the influence of CP on muscle mechanics. The passive force-stretch relationship is shown in (A) for a block of muscle tissue (C). The sarcomere length was shifted by a factor a $c_{\text{sarco}}$ and the volume fraction of the ECM was varied to represent the effects of CP. Dotted lines represent the curves corresponding to a CP base material. The active force-stretch relationship is shown in (B) for a pennate muscle geometry (D). Here, we show the influence of CP muscle geometries on the total force produced by muscle. The dotted lines correspond to CP muscle material properties ($\alpha = 0.4$, $\beta = 0.2$, and $c_{\text{sarco}} = 0.0$), while the solid lines correspond to TD muscle properties ($\alpha = 0.02$, $\beta = 0.1$, and $c_{\text{sarco}}= 0.0$)}
    \label{fig:CPExperiment}
\end{figure}

We see that our model behaves as expected for typically developed  in both the passive and active setting \autoref{fig:CPExperiment}. To investigate the passive mechanics of CP, we look at the effects of varying $c_{\text{sarco}}$ and $\alpha$ independently. We note this  cannot be done {\it in vivo}, as both of these properties vary simultaneously during CP \cite{ref:smith2011}. Using our model we are able to determine that there is a larger effect from the ECM relative to the sarcomere length in the passive setting (\autoref{fig:CPExperiment}). The benefit of our three dimensional continuum model is that we can customize the sample geometry to more accurately capture the deformation of muscle. In \autoref{fig:CPExperiment} (B), we utilize a pennate muscle geometry (\autoref{fig:CPExperiment} D), where the muscle fibres are orientated on an angle and connected to a thin stiff sheet of aponeurosis material. Using this geometry allows us to investigate to common architectural changes in CP muscle, such the decrease in physiological cross-sectional area  and volume \cite{Barrett2010}, on the muscle mechanics. We ask the question: which of the changes  \--- those to the physiological cross-section, or to the properties of the skeletal muscle itself \--- impact the force output more? To answer this, we independently varied the physiological cross section and the tissue properties. In the computational experiments corresponding to typical muscle (TD) there is an  $\approx 30\%$ decrease in force when the physiological cross section is decreased by 30\%.(\autoref{fig:CPExperiment}). This is to be expected, since it is often assumed that physiological cross section is proportional to muscle force \cite{Lieber2000}. When looking at the CP muscle tissue we instead find that there is less effect from the physiological cross section ($\approx 22\%$) (\autoref{fig:CPExperiment}), which indicates an influence from the CP material properties. Our model is able to capture the influence of changes to the material properties and architecture of CP muscle, and determine their influence on the overall muscle mechanics. Due to the highly complex structure and variability that occurs in CP affected muscle, this can not be done using experimental techniques.

\subsection{MRI-derived geometry}\label{sec:mir-geometry}
\begin{figure}[!ht]
	\centering\includegraphics[height=0.4\textheight, width=0.7\textwidth]{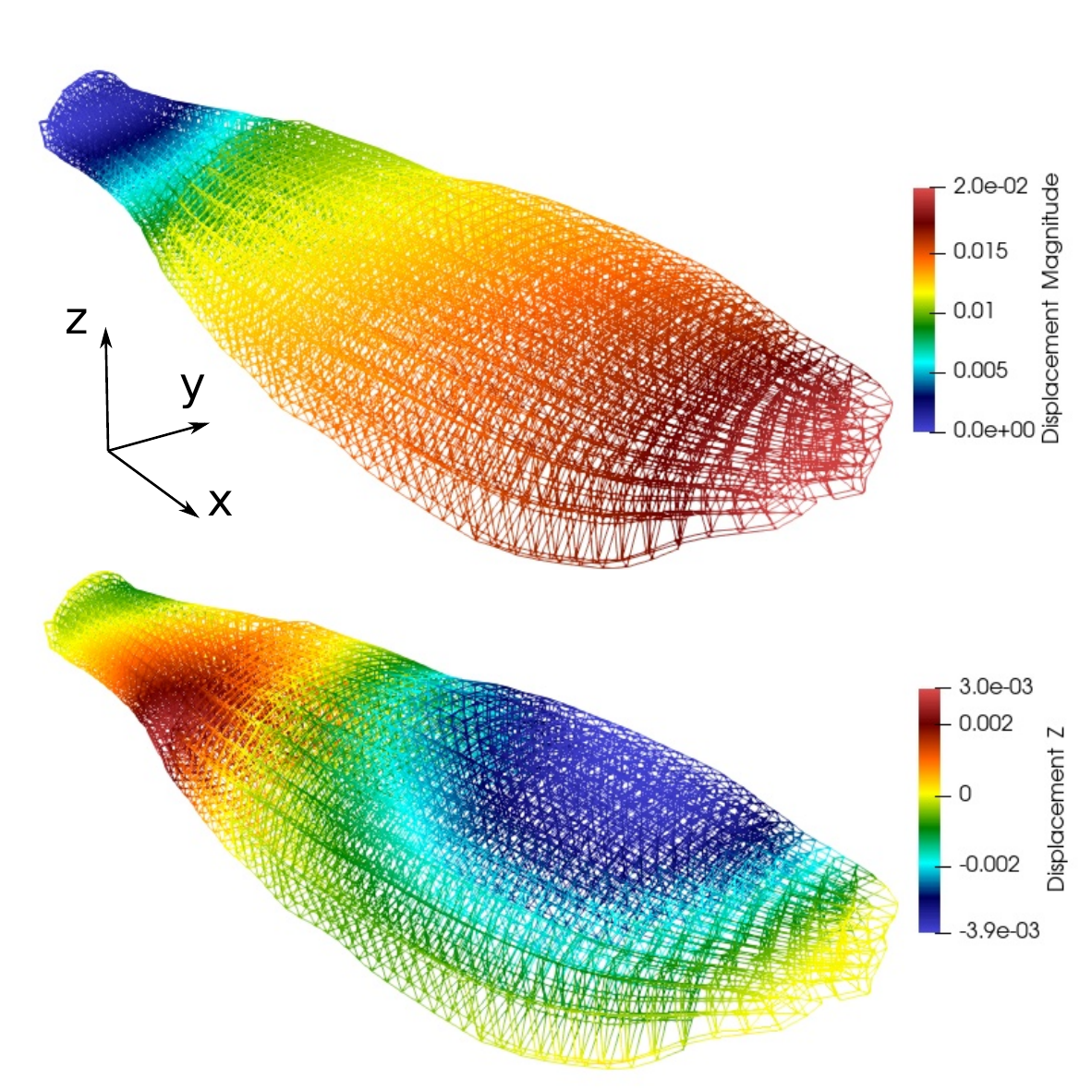}
	\caption{Magnitude of displacement (Top) and $z$ direction displacement at 100\% activation on a subject specific MRI-derived geometry. Displacements shown are in units of meters.}\label{fig:MRI-mesh}
\end{figure}

We can utilize a realistic MRI-derived geometry with our model to solve our model for subject-specific data (\autoref{fig:MRI-mesh}), which were obtained by D'Souza et al. \cite{DSouza2019}. Using these data allows us to account for variations in the structure of muscle that may vary between subjects. This could include changes to muscle thickness, length, and volume, all of which will have an effect on the mechanics of skeletal muscle \cite{ref:lieber2000}. \autoref{fig:MRI-mesh} shows a simulation where the muscle was stretched in the $x$ direction to a whole muscle stretch of 1.1, and subsequently activated to 100\% activation at a fixed length. The $+z$ and $-z$ faces have a portion of the surface covered by an aponeurosis, which the muscle fibres insert into. The fibres in this model run at 25 degrees from the $x$ axis.

In order to use our computational framework with this more realistic geometry we need the information from post-processed MRI scans to be able to define the mesh of such geometries. We also need DTI scans to be able to map the muscle fibre orientation within such a geometry; for more details see, for instance \cite{ref:wakeling2020,Bolsterlee2021}. As one can imagine, the fibre directions inside a muscle can be quite complex. Fortunately, we only need a snapshot of the geometry and the fibre orientation vectors to be used at the beginning of the simulation, as part of the initial setup of the system.

\subsection{Validation of models}
A common (and legitimate) criticism expressed by scientists in physiology of complex mathematical models concerns their validity: if many parameters are involved in the model which are fit to experimental data, how do we guard against over-fitting? In our framework, we have fit parameters describing tissue properties. We can then compare quantities such as displacement, fibre orientations, energies, and forces which are {\it not} fit; these comparisons are with {in situ} experiments. Naturally, our simulations cannot reproduce the exact values seen in experiment since we necessarily work with idealized geometries and configurations. We can, however, compare qualitative trends between experimental data and computational results.

As a first such validation, in \cite{ref:wakeling2020} we compared MRI data on muscle bulging during fixed-end contractions, with a simulation based on our model on an MRI-derived mesh. Notably, even though the location of the aponeurosis in our computation was only indirectly inferred from imaging data, and though our computations did not include the effect of tendons or adjacent muscles, our simulations predicted magnitudes of both inward and outward bulging which were consistent with those experimentally measured (Figure 12 in \cite{ref:wakeling2020}.The data we used to obtain experimental parameters did not include muscle bulging information.

We also found similar patterns of anisotropy in transverse bulging as previous ultrasound measures \cite{ref:randhawa2018}. Further, an ultrasound study of human ankle plantarflexors showed that muscles with lower initial fibre pennation angles (less than 15 degrees) tend to bulge, whereas more pennate muscles tend to thin as they activate \cite{ref:randhawa2013}, a feature of muscle contraction that we were able to replicate with our model in \cite{ref:wakeling2020}.

Next, in \cite{ref:ryan2020}, we considered fixed-end contractions in a parallel-fibered muscle which was subject to a compressive force transverse to the fibers, and the fibres are then activated. Our results showed that before activation, compression would cause muscle thickness to decrease and the fibres to change their pennation angle; this matched the trends reported in \cite{ref:ryan2019}. Our simulations also revealed a complex dependance of the longtitudinal force (oriented between the fixed ends), the length of the muscle block, and the initial fibre pennation angle. The trend we observed: compressive transverse loading lead to an increased longtitudinal force at longer lengths and a reduction at shorter lengths, consistent with experiments reported by \cite{sleboda2020}. 

Finally, to investigate the importance of mass and muscle size on the dynamics, in situ experiments were conducted on the rat plantaris muscle, in which the muscle was subjected to cyclic contractions, coupled with neuronal excitations to cause activation in synch \cite{ross2020}. Accelerations in the $x$ direction, $a_{mid}$ and $a_{end}$ were measured at the midpoint $L_{mus}/2$ and at the end point of the muscle $L_{mus}$, respectively. If the muscle mechanics of the rat plantaris muscle were genuinely well-approximated by a single fibre quasi-static Hill-type model, the scaled quantities $2a_{mid}/(L_{mus})$ and $a_{end}/L_{mus}$ should be the same. Experimentally, these are not; the difference is more pronounced when mass was added to the muscle, and also in larger muscles.

We ran simulations with the fully dynamic model on a muscle-aponeurosis configuration as in \autoref{fig:geometries}, and subjected the $+x$-face of the top aponeurosis to cyclic contractions; simultaneously, the muscle fibres were activated in a synchronized manner. We found comparable patterns of tissue accelerations in the longitudinal direction across the dynamic model during cyclic contractions \cite{ref:Ross2021} to that of the experiments on rat plantaris muscle \cite{ross2020}, again scaling the accelerations by a length scale corresponding to the location in the muscle.  We found a larger difference in scaled acceleration in the middle of the muscle relative to the end with greater muscle size for the model and greater added mass for the in situ muscle.

%% file: Manuscript/conclusions.tex
\section{Discussion and validation}\label{section:3dmuscleconclusions}
In this paper we proposed a mathematical model to describe the three dimensional deformations of skeletal muscles. We used the principles of continuum mechanics to derive this model. Muscle is considered as a nearly incompressible, fibre-reinforced and transversely isotropic composite material; internal chemical reactions are triggered to produce forces which deform the tissues without any external work done on the tissues, while passive elastic properties in the muscle fibres resist the deformation in the along fibre direction. The idea of a base material was used to encapsulate the contribution from several other tissues and fluids surrounding the fibres: blood, water, extra- and intra-cellular matrix forms part of this material. Experimental data was used to fit parameters in the model.

We also proposed a numerical approach to approximate the deformations of the tissues. A semi-implicit scheme is used to discretize the time variables whereas a conforming finite element method is used to discretize the spatial variables. As a result of such time stepping, a CFL condition governs both the time step size and mesh size in order for the computation to be stable. 
Note that the implementation of our method was not optimized to match experimental data. Instead, the outputs of our model are a result of the assumptions of the physics and mechanics in our model and its implementation.

Finally, we remark that many natural processes are not considered in this study. For example, muscle contraction in living beings produces heat. Also, fluids as water and blood form part of the overall composition of muscle tissues. Even though the influence of fluids is implicitly added to our model through the experimental data we used to fit the model parameters, the current model does not explicitly consider the effects of internal fluids as part of the deformation of the tissues. Finally, chemical reactions leading to the activation of muscle fibres is not included explicitly; only the action of these reactions is encoded in our model through the active stress $\hat{\sigma}_{act}$ and the form of the activation function which in turn depends on neuronal excitation, and can be modeled via a system of ODEs. These features correspond to an area of much current research which are not included in our model, and are left as future work.

%% file: Manuscript/acknowledgments.tex
\section*{Acknowledgments}
Sebasti\'an A. Dom\'inguez-Rivera thanks the financial support of PIMS-Canada through a PIMS postdoctoral fellowship. Nilima Nigam thanks the support of NSERC-Canada through the Discovery program.

%% file: Manuscript/note.tex
\paragraph{Note for the reader}
Please note that this manuscript represents a preprint only and has not been (or is in the process of being) peer-reviewed. A DOI link will be made available for this ArXiv preprint as soon as the peer-reviewed version is published online.

%% file: Manuscript/SIAMLS_Supplement.tex
%
%
%
%
%

\section{Supplemental information}

In this section we provide the parameters obtained through fitting described in \autoref{sec:ParameterFitting}. For completeness, we also record parameters already reported in the paper.

\subsection{Bulk moduli}
As described in \autoref{section:volumetricresponse}, each of the muscle, aponeurosis, tendon and fat components of an MTU are nearly incompressible. Aponeurosis and tendon are assumed to have the same bulk modulus, $\kappa_{apo} = \kappa_{ten} = 1\times 10^8$Pa,  \cite{ref:rahemi2014}. The bulk modulus for extracellular material is $\kappa_{ecm} = 1\times 10^6$Pa \cite{ref:blemker2005,rahemi2014regionalizing}, while that for cellular material (omitting ECM) and fat are the same $\kappa_{cell} = \kappa_{fat} = 1\times10^7 $Pa \cite{Konno2021}. The bulk modulus of muscle is modeled as a linear combination of that of fat, cellular material and ECM, depending on their respective volume fractions.

\subsection{Yeoh models of elasticity}
Each of muscle, tendon and aponeurosis are modelled as fibre-reinforced composites. The base material is isotropic and nonlinear, and the isochoric responses as modelled as in \autoref{sec:ParameterFitting}. The isochoric base materials for muscle, aponeurosis and tendon are described using Yeoh-type models; Neo-hookean models have also been tried and yield similar results \cite{ref:rahemi2014}. 
The parameters $c_k$ are obtained via non-linear regression to data from \cite{ref:mohammadkhah2016,Gillies2011decell,Alkhouli2013,Jin2013} to determine the muscle material constants and \cite{ref:mechazizi2009} for the aponeurosis and tendon base material,  \autoref{tab:BaseParam}. 
\begin{table}[!ht]
    \centering
    \begin{tabular}{|c|c|}
        \hline
         Parameter & Value/Range of Values\\\hline
         $c_{\text{1,mus}}$ & 3703  \\\hline
         $c_{\text{2,mus}}$ & -707.7 \\\hline
         $c_{\text{3,mus}}$ & 123.2  \\ \hline
         $c_{\text{1,apo}}$ & 4.6896264  \\\hline
         $c_{\text{2,apo}}$ & -3.455141 \\\hline
         $c_{\text{3,apo}}$ & 484.92055  \\ \hline
         $c_{\text{1,ten}}$ & 4.6896264  \\\hline
         $c_{\text{2,ten}}$ & -3.455141 \\\hline
         $c_{\text{3,ten}}$ & 484.92055  \\ \hline
    \end{tabular}
    \caption{List of the material parameters used in the Yeoh model for the base material for each of the tissues used in the model.}
    \label{tab:BaseParam}
\end{table}

Based on fits to experimental measurements on human breast tissue \cite{tanner2006}, the strain energy of fat is modelled as 
$ \Psi_{fat}(\bar{I_1})= 0.13\times 10^6(\bar{I}_1-3)$ Pa, \cite{ref:rahemi2015}.

\subsection{Along-fibre properties}
As discussed in \autoref{sec:ParameterFitting}, we initially fit  passive, active force-length and active force-velocity along-fibre properties of muscle to experimental data using Bezier splines, \cite{ref:ross2018}. These have directly been used in 1-D dynamic simulations. For the fully 3-dimensional calculations, these expressions are less convenient. The along-fibre passive response of muscle is modelled via a piecewise cubic polynomial which is a least-squares best fit to the Bezier curve and which has smooth derivatives:
\begin{equation}
  \hat\sigma_{\text{pass}}(\lambda) =
  \left\{
  \begin{aligned}
  & 0 && 0\leq \lambda\leq1.0\\
  &2.353(\lambda-1.0)^2 && 1.0\leq\lambda\leq1.25\\
  &3.44(\lambda-1.25)^2 + 1.18(\lambda-1.25) + 0.147 &&1.25\leq\lambda\leq1.5\\
  &0.427(\lambda-1.5)^2 + 2.90(\lambda-1.5) + 0.656 &&1.5\leq\lambda\leq1.65\\
  & 3.023(\lambda-1.65) + 1.1 &&\lambda>1.65,
  \end{aligned}\right.
\end{equation}

The active force-length response for muscle fibres is given by piecewise trigonometric fit to the Bezier splines in \cite{ref:ross2018}.
\begin{equation}
\hat\sigma_{\text{act}}(\lambda) =
\left\{
\begin{aligned}
&(0.642\sin(1.29 \lambda + 0.629)\\
                      & +0.325 \sin(5.31 \lambda -4.52)\\
                      & +0.328 \sin(6.74 \lambda + 1.69)\\
                      & +0.015 \sin(19.8 \lambda -7.39) &&\text{if }0.4\leq\lambda\leq1.75\\
                      & +0.139 \sin(8.04 \lambda + 2.54)\\
                      & +0.0018 \sin(32.2 \lambda -6.45)\\
                      & +0.012 \sin(23.2 \lambda -2.64))\\
                    & && \\
0&  &&\text{otherwise}.
\end{aligned}\right.
\end{equation}
Both the passive response and the active force-length response are based on fits to experimental data on rabbit hindlimb muscles, \cite{ref:winters2011}.

The active force-velocity relationship for muscle fibres is given by 
\begin{equation}
  \hat\sigma_{\text{vel}}(\bar\epsilon) =
  \left\{
  \begin{aligned}
  & 0 && 0\leq \bar\epsilon\leq x_1\\
  & 0.2579(\bar\epsilon-x_1)^3 + 0.1431(\bar\epsilon-x_1)^2 && x_1\leq\bar\epsilon\leq x_2\\
  & 29.8255(\bar\epsilon-x_2)^3 -0.9435(\bar\epsilon-x_2)^2 + 0.9703(\bar\epsilon-x_2)^1 + 0.3503 && x_2\leq\bar\epsilon\leq x_3\\
  &-3165.6847(\bar\epsilon-x_3)^3 + 186.1961(\bar\epsilon-x_3)^2 + 6.0908(\bar\epsilon-x_3)^1 + 1 && x_3\leq\bar\epsilon\leq x_4\\
  &0.6882(\bar\epsilon-x_4)^3 -1.4139(\bar\epsilon-x_4)^2 + 0.9678(\bar\epsilon-x_4)^1 + 1.3743 && x_4\leq\bar\epsilon\leq x_5\\
  & 1.5950 && \bar\epsilon \geq x_5,
  \end{aligned}\right.
\end{equation}
with $x_1 = -1.2 \bar\epsilon_0$, $x_2 = -0.25 \bar\epsilon_0$,$ x_3 = 0$, $x_4 = 0.05 \bar\epsilon_0$, and $x_5 = 0.75 \bar\epsilon_0$, and this is fit to experimental data on the rat flexor hallucis brevis muscle in \cite{Roots2007}.

Fibres in the aponeurosis and tendon cannot be activated. The passive along-fibre properties for the aponeurosis and tendon are the same, and are given by a cubic polynomial fit to experimental data on the human gastrocnemius reported in \cite{magnussen2003}, see also \cite{ref:rahemi2014}.
\begin{equation}
\hat\sigma_{\text{pass,apo/ten}}(\lambda) =
\left\{
\begin{aligned}
& 0 && 0\leq \lambda\leq1.0\\
&515.882034(\lambda-1.0)^2 + 0.01(\lambda - 1.0) + 0.01 &&1.0\leq\lambda\leq1.01\\
&600.590242(\lambda-1.01)^2 + 0.327640(\lambda-1.01) + 0.06168820 &&1.01\leq\lambda\leq1.02\\
&-9.975321(\lambda-1.02)^2 + 22.3394455(\lambda-1.02) + 0.2250236 &&1.02\leq\lambda\leq1.15\\
& 19.7458618(\lambda-1.15) + 2.960568 &&\lambda>1.15,
\end{aligned}\right.
\end{equation}